\documentclass[11pt]{article}

\usepackage{mathrsfs, multicol, xfrac, amsmath, amsfonts, amssymb, tikz, graphicx, wrapfig, xcolor, amsthm, epstopdf, comment}
\usepackage[caption=false]{subfig}
\usetikzlibrary{patterns}
\usepackage[margin=1in, nohead]{geometry}

\usepackage{enumitem}
\makeatletter
\newcommand{\mylabel}[2]{#2\def\@currentlabel{#2}\label{#1}}
\makeatother

\usepackage{pgfplots}
\pgfplotsset{compat=1.4}

\usepackage[colorlinks=true,breaklinks=true,bookmarks=true,urlcolor=blue,
citecolor=blue,linkcolor=blue,bookmarksopen=false,draft=false]{hyperref}

\usepackage{thmtools}
\usepackage{thm-restate}

\usepackage{tikz}
\usepackage{bbm}
\usepackage[english]{babel}
\usepackage{hyperref}
\usepackage{multirow}
\usepackage{setspace}
\usepackage[colorinlistoftodos]{todonotes}
\usepackage{array}
\usepackage[title]{appendix}
\usepackage{authblk}

\newcommand{\ubar}[1]{\mkern2mu\underline{\mkern-2mu #1\mkern-1.5mu}\mkern2mu}
\DeclareMathOperator{\preceqq}{\ubar{\preceq}}
\DeclareMathOperator{\npreceqq}{\ubar{\not\preceq}}
\DeclareMathOperator{\suceqq}{\ubar{\succeq}}
\DeclareMathOperator{\Max}{\mathrm{Max}}
\DeclareMathOperator{\Min}{\mathrm{Min}}

\newcommand{\reals}{\mathbb{R}}
\newcommand{\Z}{\mathbb{Z}}
\newcommand{\mc}{\mathcal}

\newcommand{\edittwo}[1]{{\color{black} #1}}
\newcommand{\edit}[1]{{\color{black}#1}}

\newtheorem{theorem}{Theorem}

\newtheorem{proposition}{Proposition}
\newtheorem{remark}{Remark}
\newtheorem{lemma}{Lemma}
\newtheorem{corollary}{Corollary}
\newtheorem{definition}{Definition}
\newcounter{example}[section]
\newenvironment{example}[1][]{\refstepcounter{example}\par\medskip
   \noindent \textbf{Example~\theexample. #1} \rmfamily}{\medskip}

\begin{document}

\title{Relaxations and Duality for Multiobjective Integer Programming}

\author[1]{\small Alex Dunbar}
\author[2]{Saumya Sinha}
\author[2]{Andrew J. Schaefer}
\affil[1]{\small Department of Mathematics, Emory University, Atlanta, GA}
\affil[2]{ Department of Computational Applied Mathematics and Operations Research, Rice University, Houston, TX}

\date{}

\maketitle

\begin{abstract}
Multiobjective integer programs (MOIPs) simultaneously optimize multiple objective functions over a set of linear constraints and integer variables. 
In this paper, we present continuous, convex hull and Lagrangian relaxations for MOIPs and examine the relationship among them. The convex hull relaxation is tight at supported solutions, i.e., those that can be derived via a weighted-sum scalarization of the MOIP. At unsupported solutions, the convex hull relaxation is not tight and a Lagrangian relaxation may provide a tighter bound. Using the Lagrangian relaxation, we define a Lagrangian dual of an MOIP that satisfies weak duality and is strong at supported solutions under certain conditions on the primal feasible region. \edittwo{We include a numerical experiment to illustrate that bound sets obtained via Lagrangian duality may yield tighter bounds than those from a convex hull relaxation.} Subsequently, we generalize the integer programming value function to MOIPs and use its properties to motivate a set-valued superadditive dual that is strong at supported solutions. We also define a simpler vector-valued superadditive dual that exhibits weak duality but is strongly dual if and only if the primal has a unique nondominated point.

\end{abstract}


\section{Introduction} \label{sec:intro}
Multiobjective optimization problems optimize multiple objective functions over a common set of constraints. When the constraints and objective functions are linear and the variables are constrained to be integers, the resulting problem is called a multiobjective integer (linear) program (MOIP); see \cite{Ehrgott2005MultiCriteria} for a detailed introduction. In this paper, we present relaxations and dual formulations for MOIPs. 

Duality is an extensively studied area of optimization that is used in both theoretical analysis and the development of solution methods. 
The fundamental idea in duality is to formulate a related (dual) optimization problem that can provide bounds on the original (primal) problem. For single-objective integer programs (IPs), several notions of duality have been proposed in the literature \cite{Fisher1981Lagrange,Geoffrion1974Lagrange,Hooker2009,wolsey1981integer,wolsey2014integer}. Gale, Kuhn, and Tucker first studied multiobjective linear programming (MOLP) duality in the 1950s \cite{Gale1951LinearProgrammingGames}, and the subject continued to receive attention in subsequent decades
\cite{CORLEY1984MATRIXDUAL,GurionLuc2014MOLPLagrange,Hamel2004MOLPLagrange,Heyde2008MOLPGeoDual,Heyde2009SetDual,isermann1978some,kornbluth1974duality,Rodder1977GeneralizedSaddlepoint}. Recent work on MOIPs has explored bound sets that use relaxations or variations of the original problem to bound nondominated points in the objective space \cite{BOLAND2017OptimizeoverEfficient,CERQUEUS2015SurrogateBoundSets,ehrgott2001bounds,EHRGOTT2007BoundSets,lust2010twophasePareto,machuca2016lower}. Klamroth et al. \cite{Klamroth_dual} use IP duality to derive bounds for the multiobjective problem. To our knowledge, however, a multiobjective duality framework for MOIPs has not been explored in the literature. In this paper, we extend Lagrangian and superadditive duality from single-objective IPs to the multiobjective case. 

The paper is organized as follows: in Section \ref{sec:background}, we formally define an MOIP and some related concepts of vector and set comparison. We also provide a brief review of duality for IPs and MOLPs, and an overview of bound sets. In Section \ref{sec:relaxation}, we present the Lagrangian relaxation of an MOIP and compare it with the continuous and convex hull relaxations. In Section \ref{sec:Lagdual}, we extend Lagrangian duality to the multiobjective case \edittwo{and illustrate the performance of Lagrangian dual bound sets via a numerical experiment}. In Section \ref{sec:superadditive}, we develop superadditive duals for MOIPs. We present concluding remarks in Section \ref{sec:conclusion}. 

In single-objective optimization problems, the objective function is a scalar-valued function whose maximum over the feasible region is easily defined. In contrast, a $k$-objective optimization problem has a vector-valued objective function that maps points in the feasible region to vectors in $\reals^k$. In this case, the notion of maximization is ambiguous as there may be several objective values that are mutually incomparable. Therefore, we employ the notion of Pareto efficiency and nondominance \cite{Benson2009,Ehrgott2005MultiCriteria,luc_2016_Book,Teghem2009}. We first introduce some notation for the usual element-wise order on $\reals^k$.

\begin{definition} For $x,y \in \reals^k$, define
\begin{enumerate}[label=(\roman*)]
 \item $x < y$ if $y-x$ has all positive components,
 \item $x \leqq y$ if $y-x$ has all nonnegative components, and
 \item $x \leq y$ if $x\leqq y$ but $x \neq y$.
\end{enumerate}
\end{definition}

\begin{definition}[Nondominance]
\label{def:pareto}
Given a set $S \subseteq \reals^k$, $s \in S$ is said to be nondominated (from above) if there does not exist $t \in S$ such that $s \leq t$.
\end{definition}

A set may have multiple nondominated elements that are mutually incomparable, and we denote the set of all such elements by $\Max(S)$. The set $\Min(S)$ is analogously defined by reversing the inequality in Definition \ref{def:pareto}.
%
%
Under this definition, solving a multiobjective maximization problem amounts to finding the nondominated points of the set of feasible objective function values. 

An MOIP is defined via linear constraints and objective functions.
Let $A \in \mathbb{R}^{m \times n}$ be the constraint matrix with right-hand-side $b \in \mathbb{R}^{m}$. Let $C \in \mathbb{R}^{k \times n}$ be the cost-matrix whose $i$-th row comprises the coefficients
of the $i$-th (linear) objective function. An MOIP with these $k$ objectives is defined as 
\begin{equation} 
\label{eq:MOIP}\tag{MOIP} 
\begin{aligned}
\max \enskip & Cx\\
\text{s.t.} \enskip & Ax \leqq b,\\
&x \in \mathbb{Z}^{n}_\geqq,
\end{aligned}
\end{equation}
where $\mathbb{Z}^n_{\geqq}$ is the set of integer vectors with nonnegative components. The set of positive integer vectors is denoted by $\mathbb{Z}^n_{>}$; $\reals^n_\geqq$ and $\reals^n_{>}$ are defined similarly. 
Let $\mathcal{X}$ denote the feasible region of \eqref{eq:MOIP}, 
and $\mathcal{Y} = \{Cx \ \vert \ x \in \mathcal{X} \}$ be the set of feasible objective values. The set $\mathcal{X}$ is an integer polyhedron just as in single-objective IP, while $\mathcal{Y}$ is a subset of $\reals^k$.
Solving \eqref{eq:MOIP} amounts to finding the set of nondominated points $\Max(\mathcal{Y})$, known as its {\it nondominated set}. Feasible solutions corresponding to these points are called {\it efficient solutions}.

\begin{definition}[Supported Solution]
An efficient solution $x^*$ is said to be a supported solution if there exists a scalarizing vector $\mu \in \reals^k_{>}$ such that $x^*$ is an optimal solution to the scalarized problem $\max\ \{\mu^\top Cx\ \vert \allowbreak \ Ax \leqq b, x \in \Z^n_\geqq\}$.
\end{definition}

\begin{remark}
\label{rem:weightedsum}
Henceforth, we use the term ``scalarization'' to refer to the {\it weighted-sum} scalarization $\max\ \{\mu^\top Cx\ \vert \ Ax \leqq b, x \in \Z^n_\geqq\}$, where $\mu \in \reals^k_{>}$.
\end{remark}

\begin{lemma}[\cite{Iserman74ProperEfficiency}]
\label{lem:MOLP_scal_iff}
A point $x^*$ is an efficient solution to the MOLP $\max\ \{ Cx \ \vert \allowbreak \ Ax \leqq b, x \in \reals^n_\geqq \}$ if and only if there is a vector $\mu \in \reals^k_{>}$ such that $x^*$ is an optimal solution to the scalarized problem $\max\ \{ \mu^\top Cx\ \vert \ Ax \leqq b, x \in \reals^n_\geqq\}$.
\end{lemma}

Geoffrion \cite{GEOFFRION1968ProperEfficiency} proves a general statement about scalarization for multiobjective concave functions over a convex set. Lemma \ref{lem:MOLP_scal_iff} applies this to MOLPs, asserting that finding efficient solutions to an MOLP is equivalent to finding optimal solutions to its scalarizations. This, however, is not true for the integer case -- while optimal solutions to the scalarized MOIP are efficient for the original problem, not all efficient solutions of the MOIP can be recovered in this manner \cite{Ehrgott_Scalarization,Teghem2009}.

\subsection{Extended Power Set and Set Comparison} \label{sec:setorder}
When the set of objective values $\mathcal{Y}$ is nonempty and bounded above, $\Max(\mc{Y})$ is a well-defined subset of $\reals^k$. However, this definition is ambiguous if $\mathcal{Y}$ is empty or unbounded above.
To distinguish between the two cases, we define an extended power set of $\reals^k$ (analogous to the extended real line) that contains two additional ``sets'' $\pm M_{\infty}$. 

Define a set $S \subseteq \reals^k$ to be bounded above if there exists $z \in \reals^k$ such that $s \leq z$ for all $s \in S$.
Suppose $\mathcal{Y}$ is nonempty and not bounded above. Then, we say that the maximization problem is unbounded and denote $\Max(\mc{Y}) = M_{\infty}$.
On the other hand, if $\mathcal{Y} = \emptyset$, 
which corresponds to the problem being infeasible, 
we define $\Max(\mc{Y}) = -M_{\infty}$. 
The roles of $\pm M_{\infty}$ are reversed in case of minimization.
Thus, $\pm M_{\infty}$ are analogous to $\pm \infty$ in single-objective optimization. 
These definitions extend naturally to any class of multiobjective optimization problems whose sets of feasible objective values are closed.

Unlike scalar optimization problems whose optimal values can be directly compared, nondominated sets for multiobjective problems are collections of vectors in $\reals^k$. 
Several authors \cite{CERQUEUS2015SurrogateBoundSets,EHRGOTT2007BoundSets,Lohne2011,PRZYBYLSKI2017MOBranchBound} have proposed set-orderings to compare them, especially in the context of bound sets (see Section \ref{sec:MODualBoundSets}). We employ the following set-ordering proposed by Ehrgott and Gandibleux \cite{ehrgott2001bounds}.

\begin{definition}[Pareto Set-Ordering \cite{ehrgott2001bounds}] \label{def:Set_Order}
Let $S,T \subseteq \mathbb{R}^k$ be nonempty. Define $S \preceqq T$ if 
\begin{itemize}
 \item[(i)] for every $s \in S$, there exists $t \in T$ such that $s \leqq t$, and 
 \item[(ii)] for every $s \in S$ and $t \in T$, $t \not \leq s$.
\end{itemize}
\end{definition}

In other words, every element of $S$ is dominated by an element of $T$, but no element of $S$ dominates an element of $T$ (except possibly itself). 
This ordering extends to $\pm M_{\infty}$ (see Appendix \ref{sec:app_setorder} for details). 
We then consider the following subset of the extended power set of $\reals^k$.
%
%
\begin{align}
\label{eq:setE}
\mathcal{E} = \{S \subseteq \mathbb{R}^k \ \vert\ S \not = \emptyset, \text{ if } s,t \in S \text{, then } s \not \leq t \text{ and } t \not \leq s\} \cup \{-M_{\infty},M_{\infty}\}.
\end{align}

When $k = 1$, $\mathcal{E}$ is the collection of singleton subsets of the extended real line, and the set-ordering ``$\preceqq$'' coincides with the standard ordering on $[-\infty,\infty]$. 
Proposition \ref{prop:efficient_face_in_E} implies that nondominated sets for the multiobjective optimization problems we consider in this paper will always belong to the family $\mathcal{E}$.
Proposition \ref{prop:set_partial_order} asserts that the relation ``$\preceqq$'' defines a partial order on $\mathcal{E}$. Together, these results will allow us to use this set-ordering to compare the nondominated sets of multiobjective optimization problems. 
Proofs of results from this subsection are provided in Appendix \ref{sec:app_setorder}.
\begin{restatable}{proposition}{PropMaxSinE}
\label{prop:efficient_face_in_E}
Let $S \subseteq \mathbb{R}^k$ be nonempty. If $S$ has points that are nondominated from above, then $\Max(S) \in \mathcal{E}$. If $S$ has points that are nondominated from below, then $\Min(S) \in \mathcal{E}$.
\end{restatable}

\begin{restatable}{proposition}{PropSetPartialOrder}\label{prop:set_partial_order}
The relation $\preceqq$ defines a partial order on $\mathcal{E}$. 
\end{restatable}
\noindent
The following properties of the set-ordering will be used in later proofs.
\begin{restatable}{lemma}{LemOrderSubset}
\label{lem:order_subset}
Let $T \subseteq S \subseteq \reals^k$ be nonempty sets, and let $U \in \mathcal{E}$ such that $S \preceqq U$. Then, $T \preceqq U$.
\end{restatable}

\begin{restatable}{lemma}{LemOrderMax}
\label{lem:order_max}
Let $S \subseteq \mathbb{R}^k$ be a closed set and let $U \in \mathcal{E}$ such that $S \preceqq U$. Then, $\Max(S) \preceqq U$.
\end{restatable}

\begin{remark}\label{rem:set_order_fail_min}
The analog of Lemma \ref{lem:order_max} for minimization does not hold. That is, $S \preceqq T$ does not imply $S \preceqq \Min(T)$. This is because the elements of $S$ and $\Min(T)$ may be incomparable (though we do have $s \not \geq t$ for all $s \in S$ and $t \in \Min (T)$).
For example, consider $S = \{(0,0)^\top \}$ and $T = \{(-1,1)^\top, (2,2)^\top, \allowbreak
(1,-1)^\top \}$. Then, $S \preceqq T$, but $\Min(T) = \{(-1,1)^\top, (1,-1)^\top\}$, and $S \npreceqq \Min(T)$. 
\end{remark}

\subsection{Duality in Integer Programming}
\label{sec:IPduality}
Duality for single-objective IPs has been well studied, and we briefly review Lagrangian and superadditive duality; see \cite{wolsey2014integer} for a comprehensive discussion including the following results. Consider the IP
\begin{equation} \label{eq:IP}\tag{IP}
\begin{aligned}
z_{IP} = \max \enskip cx \quad
\text{s.t.} \enskip Ax \leqq b, \ x \in \mathbb{Z}^n_\geqq.
\end{aligned}
\end{equation}

Suppose the constraint matrix $A$ is composed of two sub-matrices, $A^1 \in \mathbb{R}^{m_1 \times n}$ and $A^2 \in \mathbb{R}^{(m-m_1)\times n}$, so that the constraints corresponding to $A^1$ are the so-called complicating constraints, while those corresponding to $A^2$ are somewhat easier to handle. Then, a Lagrangian dual is constructed by penalizing the complicating constraints by means of a dual variable $\lambda$ as follows.
\begin{equation}
\label{eq:IP_Dual}
 z_{\mathrm{LD}} = \min_{\lambda \in \mathbb{R}^m_{\geqq}}\ \max_{x \in Q}\ (cx + \lambda^\top(b^1 - A^1x)),
\end{equation}
\noindent 
where $Q = \{x \in \mathbb{Z}_\geqq^n\ \vert \ A^2x \leqq b^2\}$, and $b^1$ and $b^2$ are the sub-vectors of $b$ corresponding to the rows included in $A^1$ and $A^2$, respectively. The Lagrangian dual \eqref{eq:IP_Dual} satisfies the following properties \cite{Fisher1981Lagrange,Geoffrion1974Lagrange}.

\begin{proposition}[Weak Lagrangian Duality]\label{prop:IP_weak_Lag_Dual}
For each $\lambda \in \reals^{m_1}_+$ and each $x$ feasible to \eqref{eq:IP}, 
\begin{equation*} 
\begin{aligned}
cx \leq cx + \lambda^\top(b^1 - A^1x).
\end{aligned} 
\end{equation*}
It follows that $z_{IP} \leq z_{\mathrm{LD}}$.
\end{proposition}

\begin{theorem}\label{thm:LD_as_LP}
The optimal value $z_{\mathrm{LD}}$ of \eqref{eq:IP_Dual} is equal to the optimal value of the following linear program (LP):
\begin{equation*} \begin{aligned}
\max \enskip cx \quad
\text{ \emph{s.t.}} \enskip A^1x \leqq b^1, \ x \in \text{\emph{conv}}(Q).
\end{aligned} \end{equation*}
\end{theorem}

\begin{theorem}\label{thm:Single_objective_LD_iff}
For fixed $A^1,b^1$ and $Q$, the Lagrangian dual \eqref{eq:IP_Dual} is strong for any cost vector $c$ if and only if 
\begin{equation*} \begin{aligned}
\text{\emph{conv}}\left(Q \cap \{x \in \mathbb{R}^n_\geqq\ \vert \ A^1x \leqq b^1\}\right) = \text{\emph{conv}}(Q)\cap\{x \in \mathbb{R}^n_\geqq\ \vert \ A^1x \leqq b^1\}.
\end{aligned} \end{equation*}
\end{theorem}

The value function of an IP is defined as its optimal value parameterized by the constraint right-hand-side $b$. That is, $
z(b) = \max\{cx\ \vert \ Ax \leqq b, x \in \mathbb{Z}^n_\geqq\}$.
The value function is nondecreasing and superadditive over its domain, which motivates the following superadditive dual for \eqref{eq:IP}. 
\begin{equation}
\label{eq:Single_Objective_Super}
\begin{aligned}
\min \enskip & F(b)\\
\text{s.t.} \enskip & F(A_j) \geqq c_j, & j=1, \ldots, n, \\
&F(0) = 0,\\
&F:\mathbb{R}^m \rightarrow \mathbb{R}, & \text{nondecreasing and superadditive}.
\end{aligned}
\end{equation}
\noindent The superadditive dual \eqref{eq:Single_Objective_Super} satisfies weak and strong duality \cite{JEROSLOW1978121,wolsey1981integer}.

\begin{proposition}[Weak Superadditive Duality]\label{prop:IP_weak_SA}
If $x$ is feasible for \eqref{eq:IP} and $F$ is feasible for \eqref{eq:Single_Objective_Super}, then $cx \leq F(b)$. If \eqref{eq:IP} is unbounded, then \eqref{eq:Single_Objective_Super} is infeasible.
\end{proposition}

\begin{theorem}[Strong Superadditive Duality]
\label{thm:single_obj_strong_super}
If \eqref{eq:IP} has an optimal solution $x^*$ with $cx^* < \infty$, then \eqref{eq:Single_Objective_Super} has an optimal solution $F^*$ with $F^*(b) = cx^*$. 
\end{theorem}

\subsection{Duality in Multiobjective Linear Programming}

We now provide a brief overview of MOLP duality; see Luc \cite{Luc2011MOLPDual} for a survey. Consider the MOLP
\begin{equation}
\label{eq:MOLP_gen}
 \begin{aligned}
 \max \enskip Cx \quad
 \text{s.t. } \enskip Ax = b, \
  x \in \reals^n_{\geqq}.
 \end{aligned}
\end{equation}

A matrix optimization dual for \eqref{eq:MOLP_gen} was first proposed by Gale et al. in \cite{Gale1951LinearProgrammingGames}. Several types of MOLP duality have since been explored, such as geometric duality \cite{Heyde2008MOLPGeoDual} and set valued duality \cite{Heyde2009SetDual}. 
Isermann \cite{isermann1978some} proposed a vector-valued dual problem over a space of matrices as an extension of the single-objective LP dual. This dual is strong in the sense that if both the primal and dual problems are feasible, then they \edit{have} nondominated points that coincide. Isermann's dual has also been extended to other variants by considering different feasible sets for the dual matrix variables, such as Corley's duality in \cite{CORLEY1984MATRIXDUAL}.


Lagrangian duality for MOLPs was introduced by Hamel et al.  \cite{Hamel2004MOLPLagrange} and extended by Gourion and Luc \cite{GurionLuc2014MOLPLagrange}. 
Given the MOLP \eqref{eq:MOLP_gen}, define a set-valued function $G$ on $\reals^{k\times m}$ as 
$G(\Lambda) = \Max\ (\{Cx + \Lambda(b-Ax)\ \vert \ x \geqq 0\})$. 
Then, the Lagrangian minmax problem $\min \{G(\Lambda)\ \vert \ \allowbreak \Lambda \in \reals^{k \times m}\}$ is expressed as
\begin{equation}\label{eq:MOLP_LD}
\begin{aligned}
 \min &\enskip \Lambda b + (C - \Lambda A)u\\
 \text{s.t.} & \enskip \alpha^\top \Lambda A \geqq \alpha^\top C,\\
 & \enskip \alpha^\top (C - \Lambda A)u = 0, \\
 &\enskip \alpha \in \reals^k_{>}, \enskip u \in \reals^n_\geqq.
\end{aligned}
\end{equation}


\begin{lemma}[Strong duality of \eqref{eq:MOLP_LD} \cite{luc_2016_Book}]
\label{lem:MOLP_LD_STRONG}
If $\mathcal{Y}_{\mathrm{M}}$ and $\mathcal{Y}_{\mathrm{L}}$ are the sets of feasible objective values of \eqref{eq:MOLP_gen} and \eqref{eq:MOLP_LD} respectively, then $\Max(\mathcal{Y}_{\mathrm{M}}) = \Min(\mathcal{Y}_{\mathrm{L}})$.
\end{lemma}

\subsection{Bound Sets for Multiobjective Integer Programming}\label{sec:MODualBoundSets} 

While MOLP duality is well-studied, duality for MOIPs remain relatively unexplored. However, the related idea of using relaxations or variations of an MOIP to derive bounds on its nondominated points has been explored through bound sets \cite{Ehrgott2005MultiCriteria,ehrgott2001bounds,EHRGOTT2007BoundSets}. 
Bound sets are helpful in computing the entire nondominated set of an MOIP and have been popular in recent work on MOIPs \cite{BOLAND2017OptimizeoverEfficient,CERQUEUS2015SurrogateBoundSets,DACHERT2017EfficientComputationSearch,Haimes1971epsilon_constraint,Klamroth_dual,lust2010twophasePareto,machuca2016lower}; see \cite{Halffmann2022Review} for a survey of algorithms for MOIPs.

A common upper bound on $\mathcal{Y}$ is given by the {\it ideal point} defined coordinatewise as $y^I_i = \max\limits_{x \in \mathcal{X}} (Cx)_i$, $i = 1, 2, \ldots, k$. 
Similarly, a common lower bound is the \emph{nadir point} defined as $y^N_i = \min \{(Cx)_i \ \vert \ x \text{ efficient for \eqref{eq:MOIP}}\}$, $i = 1,2,\ldots, k$ \cite{Ehrgott2005MultiCriteria}. 
Then, for each $y \in \Max(\mathcal{Y})$, $y^N \leqq y \leqq y^I$. 
Because $\Max(\mathcal{Y})$ is in general a set with more than one element, the notion of bounding has been extended to sets. Fix a subset $Y \subseteq \mathcal{Y}$. 
Ehrgott and Gandibleux \cite{ehrgott2001bounds} define an upper (resp. lower) bound set $U$ (resp. $L$) as a subset of $\mathbb{R}^k$ such that $Y \preceqq U$ (resp. $L \preceqq Y$). 
Note that this definition of bound sets is not unique. For example, in \cite{EHRGOTT2007BoundSets}, Ehrgott and Gandibleux use an alternative definition that imposes additional conditions on the sets $L$ and $U$. However, in this paper, we follow \cite{ehrgott2001bounds} and define bound sets through ``$\preceqq$".

In \cite{AnejaNair79BicriteriaTransport}, the authors solve a sequence of scalarized problems to obtain an upper bound set. Their approach uses the fact that the nondominated set of the convex hull of the supported nondominated points of \eqref{eq:MOIP} is an upper bound set for $\mathcal{Y}$. 
Moreover, because the scalarization of a relaxation is a relaxation of the scalarization, this approach can be extended to problems for which the supported efficient solutions are difficult to compute \cite{EHRGOTT2007BoundSets}. This technique has further been adapted to specific problems by leveraging the problem features \cite{CERQUEUS2015SurrogateBoundSets,lust2010twophasePareto,machuca2016lower}. Alternatively, bound sets may be derived via search space decomposition, which uses local information based on regions of the search space \cite{BOLAND2017OptimizeoverEfficient,DACHERT2017EfficientComputationSearch}.
Klamroth et al. \cite{Klamroth_dual} compute bound sets using single-objective duality applied to $\varepsilon$-constraint scalarizations \cite{Haimes1971epsilon_constraint} of the MOIP.

In the subsequent sections, we provide several bounds on the nondominated points of MOIPs through relaxations and multiobjective dual problems.

\section{Relaxations for Multiobjective Integer Programs} \label{sec:relaxation}
A $k$-objective maximization problem is a relaxation of \eqref{eq:MOIP} if its feasible region contains the feasible region of \eqref{eq:MOIP}, and its nondominated set is an upper bound set for $\Max(\mc{Y})$. In this section, we first review the continuous and convex hull relaxations, and then present Lagrangian relaxation for MOIPs.

\vspace{-2mm}
\subsection{The MOLP Relaxation}\label{sec:MOLPRelax}
The continuous relaxation of \eqref{eq:MOIP} is obtained by dropping the integrality constraints, which yields the following MOLP.

\begin{equation}\label{eq:MOLP}\tag{MOLP} 
\begin{aligned} 
\max \enskip Cx \quad \text{s.t.} \enskip & Ax \leqq b, \ x \in \mathbb{R}_\geqq^{n}.
\end{aligned}
\end{equation}

In recent years, the MOLP relaxation has been used in several search-based methods for solving MOIPs \cite{FORGET2022909,jozefowiez2012generic,mavrotas2005multi,VINCENT2013498}.
The feasible region of \eqref{eq:MOLP} clearly contains the feasible region of \eqref{eq:MOIP}. Proposition \ref{prop:MOLP_relax} shows that its nondominated set provides an upper bound for $\Max(\mc{Y}_{\mathrm{MOIP}})$\footnote{Because we consider multiple optimization problems in this section, we use $\mc{X}_{\mathrm{P}}$ and $\mc{Y}_{\mathrm{P}}$ to denote the feasible region and set of feasible objective values for problem (P), respectively.}, so that \eqref{eq:MOLP} is a relaxation of \eqref{eq:MOIP}; the proof is provided in Appendix \ref{sec:app_relaxation}.

\begin{restatable}{proposition}{PropMOLPRelax}
\label{prop:MOLP_relax} 
If $\mathcal{Y}_{\mathrm{MOLP}}$ is the set of feasible objective values of \eqref{eq:MOLP}, then $\Max(\mathcal{Y}_{\mathrm{MOIP}}) \preceqq \Max(\mathcal{Y}_{\mathrm{MOLP}})$.
\end{restatable}

\subsection{The Convex Hull Relaxation}\label{sec:CHRelax}
The convex hull relaxation of \eqref{eq:MOIP} is defined as
\begin{equation}
\label{eq:CH}\tag{CH}
\begin{aligned}
\max &\enskip Cx \quad
\text{s.t.} \enskip x \in \text{conv}\left(\left\{x \in \Z^n_\geqq \ \vert \ Ax \leqq b \right\}\right).
\end{aligned}
\end{equation}
\noindent 
Again, the feasible region of \eqref{eq:CH} clearly contains the feasible region of \eqref{eq:MOIP}, and Proposition \ref{prop:MOIP_leqq_CH} asserts that its nondominated set provides an upper bound for $\Max(\mathcal{Y}_{\mathrm{MOIP}})$. Thus, \eqref{eq:CH} is a relaxation of \eqref{eq:MOIP}. The convex hull relaxation has been widely used in solution algorithms for MOIPs \cite{EHRGOTT2007BoundSets,ozpeynirci2010exact,przybylski2010recursive,sourd2008multiobjective}. We include some relevant results here; proofs are available in Appendix \ref{sec:app_relaxation}.

\begin{restatable}{proposition}{PropMOIPLeqqCH}\label{prop:MOIP_leqq_CH}
If $\mathcal{Y}_{\mathrm{CH}}$ is the set of feasible objective values of \eqref{eq:CH}, then
$\Max(\mathcal{Y}_{\mathrm{MOIP}}) \preceqq \allowbreak \Max(\mathcal{Y}_{\mathrm{CH}})$.
\end{restatable}


\begin{remark}
\label{rem:CH_MOLP}
Because $\mathcal{X}_{\mathrm{CH}}\subseteq \mathcal{X}_{\mathrm{MOLP}}$, the convex hull relaxation of an MOIP is tighter than its continuous relaxation. That is, $\Max(\mathcal{Y}_{\mathrm{CH}}) \preceqq \Max(\mathcal{Y}_{\mathrm{MOLP}})$. 
\end{remark}

The inequality in Proposition \ref{prop:MOIP_leqq_CH} does not necessarily hold with equality.
In the subsequent results, we explore connections between efficient solutions of $\eqref{eq:CH}$ and $\eqref{eq:MOIP}$. 

\begin{restatable}{proposition}{PropCHMOIPIntegral} \label{prop_CH_MOIP_integral}
Let $x^*$ be an efficient solution of \eqref{eq:CH}. Then, $x^*$ is an efficient solution of \eqref{eq:MOIP} if and only if $x^*$ is integral. 
\end{restatable}

\begin{restatable}{proposition}{PropCHIntSol} \label{prop:CH_int_sol}
If \eqref{eq:CH} has efficient solutions, then it must have an integral efficient solution. 
\end{restatable}

Propositions \ref{prop_CH_MOIP_integral} and \ref{prop:CH_int_sol} imply that solving \eqref{eq:CH} returns at least one efficient solution of \eqref{eq:MOIP}. Solutions to an MOIP obtained via its convex hull relaxation are further characterized as supported solutions in Proposition \ref{prop:scal_convex}.

\begin{restatable}[\cite{Ehrgott_Scalarization,Przybylski2010TwoPhase}]{proposition}{PropScalConvex}
\label{prop:scal_convex}
Let $x^*$ be an efficient solution of \eqref{eq:MOIP}. Then, $x^*$ is a supported solution if and only if it is an efficient solution of \eqref{eq:CH}. 
\end{restatable}

In contrast to single-objective IPs, where every optimal solution to the IP is also optimal for its convex hull relaxation,
\eqref{eq:MOIP} may have efficient solutions that are not efficient for \eqref{eq:CH}. This is because the images of  unsupported solutions do not lie on the nondominated frontier of $\text{conv}(\mathcal{Y})$. 

\subsection{Lagrangian Relaxation for MOIPs}\label{sec:LagRel}

In this section, we extend Lagrangian relaxation to MOIPs. Suppose the constraint matrix $A$ is comprised of two sub-matrices $A^1 \in \reals^{m_1 \times n}$ corresponding to ``complicating" constraints, and $A^2\in \reals^{(m - m_1) \times n}$ corresponding to ``simple" constraints. 
Let $b^1$ and $b^2$ be the corresponding sub-vectors of the constraint right-hand-side $b$.
As in the single-objective case, we only dualize the complicating constraints. 
Given a matrix of multipliers $\Lambda \in \mathbb{R}_\geqq^{k\times m_1}$, the Lagrangian relaxation of \eqref{eq:MOIP} is the following multiobjective \edit{optimization} problem.

\begin{equation}
\label{eq:Lagrangian_Relaxation} \tag{LR($\Lambda$)}
\begin{aligned}
\max &\enskip Cx + \Lambda(b^1 - A^1x)\\
\text{s.t.} & \enskip x \in Q = \{x\in \mathbb{Z}^n_{\geqq}\ \vert \ A^2x \leqq b^2\}.
\end{aligned}
\end{equation}

\noindent 
If $x$ is feasible to \eqref{eq:MOIP}, then it is feasible to \eqref{eq:Lagrangian_Relaxation} as well. Therefore, \eqref{eq:Lagrangian_Relaxation} is a relaxation of \eqref{eq:MOIP} provided its nondominated set yields an upper bound for $\Max(\mathcal{Y}_{\mathrm{MOIP}})$. This is established in Proposition \ref{prop_weak_dual}. 

\begin{proposition}
\label{prop_weak_dual} 
For $\Lambda \in \reals^{k\times m_1}_\geqq$, let $\mathcal{Y}_{\mathrm{LR(\Lambda)}}$ be the set of feasible objective values of \eqref{eq:Lagrangian_Relaxation}. Then,
$\Max(\mathcal{Y}_{\mathrm{MOIP}}) \preceqq \Max(\mathcal{Y}_{\mathrm{LR(\Lambda)}})$. 
\end{proposition}

\proof 
If \eqref{eq:MOIP} is infeasible, then $\Max(\mathcal{Y}_{\mathrm{MOIP}}) = -M_{\infty}$ and the result follows trivially. 
Next, if $x \in \mathcal{X}$, then $x$ is feasible to \eqref{eq:Lagrangian_Relaxation} and $Cx \leqq Cx + \Lambda(b^1 - A^1x)$. So, if \eqref{eq:MOIP} is unbounded above, then so is \eqref{eq:Lagrangian_Relaxation} and the result holds. 
Finally, suppose \eqref{eq:MOIP} is feasible and bounded and let $y^* = Cx^* \in \Max(\mathcal{Y}_{\mathrm{MOIP}})$.
Then $x^*$ is feasible for \eqref{eq:Lagrangian_Relaxation} as well.
If \eqref{eq:Lagrangian_Relaxation} is unbounded, then the result holds. Otherwise, there exists $\hat{y} \in \Max(\mathcal{Y}_{\mathrm{LR(\Lambda)}})$ such that $\hat{y} \geqq Cx^* + \Lambda(b^1 - A^1x^*) \geqq Cx^*$.
On the other hand, suppose there is $\tilde{y}\in \Max(\mathcal{Y}_{\mathrm{LR}(\Lambda)})$ such that $\tilde{y} \leq y^*$. This implies $\tilde{y} \leq \hat{y}$, which contradicts the nondominance of $\tilde{y}$ for \eqref{eq:Lagrangian_Relaxation}. Thus, $\Max(\mathcal{Y}_{\mathrm{MOIP}}) \preceqq \Max(\mathcal{Y}_{\mathrm{LR(\Lambda)}})$.
\qed 
\endproof

\begin{remark} \label{rem:NoRelation_CH_LR}
We noted in Remark \ref{rem:CH_MOLP} that $\Max(\mathcal{Y}_{\mathrm{CH}}) \preceqq \Max(\mathcal{Y}_{\mathrm{MOLP}})$ for any MOIP. In general, there is no such relationship between $\Max(\mathcal{Y}_{\mathrm{LR(\Lambda)}})$ and either $\Max(\mathcal{Y}_{\mathrm{CH}})$ or $\Max(\mathcal{Y}_{\mathrm{MOLP}})$. This is illustrated in Example \ref{ex:LR_better_everywhere}.
\end{remark}

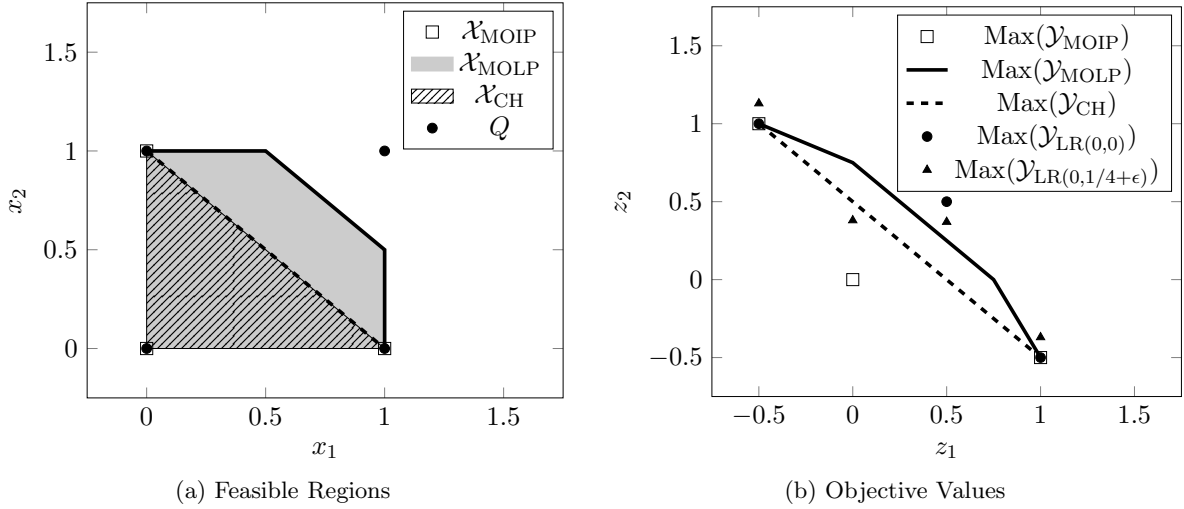
\begin{figure*}
\subfloat[Feasible Regions]{
\resizebox{0.48\textwidth}{!}{
\centering 
 \begin{tikzpicture}[scale = 0.8]
\begin{axis}[
xmin = {-0.25},
xmax = {1.75},
ymin = {-0.25},
ymax = {1.75},
xlabel = {$x_1$},
ylabel = {$x_2$}
]

\addplot[only marks, mark size = 2.5, mark = square] coordinates {(0,0) (1,0) (0,1)};
\addplot[fill, black!20!white, area legend]coordinates {(0,0) (0,1) (0.5,1) (1,0.5) (1,0)} \closedcycle;  
\addplot[pattern = north east lines, area legend] coordinates {(0,0) (0,1) (1,0)} \closedcycle;
\addplot[only marks] coordinates {(0,0) (1,1) (0,1) (1,0)};
\addplot[dashed, line width = 1.5pt] coordinates {(0,1) (1,0)};
\addplot[line width = 1.5pt] coordinates {(0,1) (0.5,1) (1, 0.5) (1,0)};

\legend{$\mathcal{X}_{\mathrm{MOIP}}$, $\mathcal{X}_{\mathrm{MOLP}}$, $\mathcal{X}_{\mathrm{CH}}$, $Q$}

\end{axis}
\end{tikzpicture} 
 \label{fig:ex_LR_better_everywhere_a}
}}
\subfloat[Objective Values]{
\resizebox{0.49\textwidth}{!}{
\centering 
\begin{tikzpicture}[scale = 0.8]
\begin{axis}[
xmin = {-0.75},
xmax = {1.75},
ymin = {-0.75},
ymax = {1.75},
xlabel = {$z_1$},
ylabel = {$z_2$}
]

\addplot[only marks, mark size = 2.5, mark = square] coordinates {(0,0) (-0.5,1) (1, -0.5)};
\addplot[line width = 1.5pt] coordinates {(-0.5,1) (0,0.75) (0.75, 0) (1, -0.5)};
\addplot[dashed, line width = 1.5pt] coordinates {(-0.5,1) (1, -0.5)};
\addplot[only marks] coordinates {(-0.5,1) (1, -0.5) (0.5,0.5)};
\addplot[only marks, mark = triangle*] coordinates {(0,0.38) (1, -0.37) (-0.5, 1.13) (0.5, 0.37)};

\legend{$\Max(\mathcal{Y}_{\mathrm{MOIP}})$,  $\Max(\mathcal{Y}_{\mathrm{MOLP}})$, $\Max(\mathcal{Y}_{\mathrm{CH}})$, $\Max(\mathcal{Y}_{\mathrm{LR}(0 \text{,} 0)})$, $\Max(\mathcal{Y}_{\mathrm{LR}(0 \text{,} 1/4 + \epsilon)})$}
\end{axis}
\end{tikzpicture}
 \label{fig:ex_LR_better_everywhere_b}
}}
\caption{Feasible regions and nondominated sets for the MOIP in Example \ref{ex:LR_better_everywhere} and its continuous, convex hull, and Lagrangian relaxations. For this example, a subset of Lagrangian upper bounds outperforms both the continuous and convex hull relaxations.}\label{fig:ex_LR_better_everywhere}
\end{figure*}

\begin{example} 
\label{ex:LR_better_everywhere}
Consider the MOIP 
\begin{equation*} \begin{aligned}
\max \enskip & \begin{bmatrix}1 & -\frac{1}{2} \\ -\frac{1}{2} & 1 \end{bmatrix} \begin{bmatrix} x_1 \\ x_2 \end{bmatrix} \quad
\text{s.t} \enskip x_1 + x_2 \leqq \frac{3}{2}, \ x_1,x_2 \in \{0,1\}.
\end{aligned} \end{equation*}

For the Lagrangian relaxations, we dualize the linear constraint. 
Then, for each $\Lambda = (\lambda_1,\lambda_2)^\top \geqq 0$, the set of feasible objective values for the Lagrangian relaxation is %
\begin{equation*} 
\begin{aligned}
\mc{Y}_{\mathrm{LR}(\Lambda)} = \left\{\begin{bmatrix}\frac{3}{2}\lambda_1 \\ \frac{3}{2}\lambda_2 \end{bmatrix}, \begin{bmatrix}1 + \frac{1}{2}\lambda_1\\-\frac{1}{2}+ \frac{1}{2}\lambda_2 \end{bmatrix}, \begin{bmatrix}-\frac{1}{2} + \frac{1}{2}\lambda_1\\1+ \frac{1}{2}\lambda_2 \end{bmatrix},\begin{bmatrix}\frac{1}{2} - \frac{1}{2}\lambda_1\\\frac{1}{2}- \frac{1}{2}\lambda_2 \end{bmatrix} \right\}.
\end{aligned} 
\end{equation*}

In particular, the points $(1,-\frac{1}{2})^\top$ and $ (-\frac{1}{2},1)^\top$ are feasible objective values for the Lagrangian relaxation with $\Lambda = (0,0)^\top$, and correspond to the supported efficient solutions $(1,0)^\top$ and $(0,1)^\top$ of the MOIP. This shows that the Lagrangian relaxation can be tight at some points in $\Max(\mc{Y}_{\mathrm{MOIP}})$.

Next, if $\Lambda = (0,\frac{1}{4} + \epsilon)^\top$ for small $\epsilon > 0$,
\[
\Max(\mathcal{Y}_{\mathrm{LR(\Lambda)}}) = \left\{\begin{bmatrix}0\\ \frac{3}{8} + \frac{3}{2}\epsilon \end{bmatrix}, \begin{bmatrix}1 \\-\frac{3}{8}+ \frac{1}{2}\epsilon \end{bmatrix}, \begin{bmatrix}-\frac{1}{2}\\\frac{9}{8}+ \frac{1}{2}\epsilon \end{bmatrix},\begin{bmatrix}\frac{1}{2}\\\frac{3}{8}- \frac{1}{2}\epsilon \end{bmatrix} \right\}.
\]

The relaxations are illustrated in Figure \ref{fig:ex_LR_better_everywhere}, which depicts $\Max(\mathcal{Y}_{\mathrm{MOIP}})$, $\Max(\mathcal{Y}_{\mathrm{CH}})$, and $\Max(\mathcal{Y}_{\mathrm{MOLP}})$, and the nondominated sets $\Max(\mathcal{Y}_{\mathrm{LR}((0,0)^\top)})$ and $\Max(\mathcal{Y}_{\mathrm{LR}((0,1/4 + \epsilon)^\top)})$ corresponding to Lagrangian relaxations. In this case, there is a subset 
\begin{equation*}
\begin{aligned}
   S & = \Big\{ \Big(-\frac{1}{2},1\Big)^\top, \Big(0,\frac{3}{8} + \frac{3\epsilon}{2}\Big)^\top, \Big(1, -\frac{1}{2}\Big)^\top \Big\} \\
   & = \Min\bigg( \Max(\mathcal{Y}_{\mathrm{LR}((0,0)^\top)}) \bigcup \Max(\mathcal{Y}_{\mathrm{LR}((0,1/4 + \epsilon)^\top)}) \bigg),
\end{aligned}
\end{equation*}

\noindent such that  $\Max(\mathcal{Y}_{\mathrm{MOIP}}) \preceqq S \preceqq \Max(\mathcal{Y}_{\mathrm{CH}}) \preceqq \Max(\mathcal{Y}_{\mathrm{MOLP}})$, and none of the inequalities holds with equality. 
\qed\end{example}

Example \ref{ex:LR_better_everywhere} demonstrates that there are problems for which Lagrangian relaxations can provide strictly tighter upper bounds on the nondominated set of \eqref{eq:MOIP} than the convex hull relaxation.

\section{Lagrangian Duality for Multiobjective Integer Programs} 
\label{sec:Lagdual}

In the previous section, we showed that for any nonnegative matrix $\Lambda$, the nondominated set of the Lagrangian relaxation provides an upper bound set for \eqref{eq:MOIP}. 
Moreover, in Example \ref{ex:LR_better_everywhere}, we found that at unsupported solutions, an appropriate choice of $\Lambda$ may yield a bound that is tighter than the convex hull relaxation. This motivates a strategy where we search for the ``best'' among all the upper bounds obtained via Lagrangian relaxation. 
To this end, we consider the set 
\begin{equation*} 
\begin{aligned}
\mathcal{Y}_{\mathrm{LD}} = \bigcup_{\Lambda \in \reals^{k\times m_1}_{\geqq}}\Max(\mathcal{Y}_{\mathrm{LR}(\Lambda)}).
\end{aligned} 
\end{equation*}
Thus, $\mathcal{Y}_{\mathrm{LD}}$ is a subset of $\reals^k$ comprised of the nondominated points of all possible Lagrangian relaxations of \eqref{eq:MOIP}. Then, a natural approach to finding the best bounds due to the Lagrangian relaxations is to consider the elements of $\mathcal{Y}_{\mathrm{LD}}$ that are nondominated from below. Therefore, we define the Lagrangian dual of \eqref{eq:MOIP} as
\begin{equation}
\label{eq:Lagrange_dual}\tag{LD}
\Min( \mathcal{Y}_{\mathrm{LD}})
= \Min \bigg( \bigcup\limits_{\Lambda \geqq 0} \Max(\mathcal{Y}_{\mathrm{LR}(\Lambda)}) \bigg).
\end{equation}
\subsection{Geometry of the Set of Lagrangian Dual Feasible Values}
\label{sec:geometry}
Before we establish properties of the Lagrangian dual \eqref{eq:Lagrange_dual}, we first analyze its set of feasible objectives $\mathcal{Y}_{\mathrm{LD}}$.
Because $\mathcal{Y}_{\mathrm{LD}}$ is an uncountable union of closed sets in $\reals^k$, it is not guaranteed to be closed (or open). This is illustrated in Example \ref{ex:P_not_open}, where we analytically derive the set $\Max(\mathcal{Y}_{\mathrm{LR}(\Lambda)})$ for all $\Lambda \geqq 0$ and use it to obtain a complete description of $\mathcal{Y}_{\mathrm{LD}}$. In general, however, an explicit description of $\mathcal{Y}_{\mathrm{LD}}$ may be difficult to derive. 

\begin{figure}
    \centering
    \resizebox{0.6\textwidth}{!}{
\begin{tikzpicture}[scale = 0.8]
\begin{axis}[
xmin = {-0.75},
xmax = {1.5},
ymin = {-0.75},
ymax = {1.25},
legend style={at={(0.5,-0.25)},anchor=north},
legend style={/tikz/every even column/.append style={column sep=0.5cm}},
legend columns = -1,
xlabel = {$z_1$},
ylabel = {$z_2$}
]

\addplot[only marks, mark size = 3, mark = square] coordinates {(0,0) (-0.5,1) (1, -0.5)};
\addplot[fill, black!20!white, area legend] coordinates{
(-0.5,0.25)
(-0.5,0.5)
(0,0.5)
(0,1.25)
(1.5,1.25)
(1.5,0)
(0.5,0)
(0.5,-0.5)
(0.25,-0.5)
(0.25,0.25)
} 
--cycle;

\addplot[dashed, line width = 1.5pt] coordinates{
(-0.5, 0.5)
(-0.5,0.25)
(0.25, 0.25)
(0.25,-0.5)
(0.5,-0.5)
};

\addplot[only marks] coordinates {
(-0.5,1) (0.25,0.25) (1, -0.5)
};

\addplot[line width = 2pt] coordinates{
(-0.5,0.5) (0,0.5) (0, 1.25)
};

\addplot[line width = 2pt] coordinates{
(0.5,-0.5) (0.5,0) (1.5, 0)
};





\legend{
{\large $\Max(\mathcal{Y}_{\mathrm{MOIP}})$},
{\large $\mathcal{Y}_{\mathrm{LD}}$},
{\large Unattained Limit Points of $\mathcal{Y}_{\mathrm{LD}}$},
{\large $\Min(\mathcal{Y}_{\mathrm{LD}})$}
}

\end{axis}
\end{tikzpicture}
}
    \caption{The set $\mathcal{Y}_{\mathrm{LD}}$ comprising the nondominated points of Lagrangian relaxations of the MOIP in Example \ref{ex:P_not_open}. $\mathcal{Y}_{\mathrm{LD}}$ is non-convex, disconnected, and neither closed nor open. Moreover, the limit point of $\mathcal{Y}_{\mathrm{LD}}$ at $(\frac{1}{4}, -\frac{1}{2})^\top $ is dominated by the primal nondominated point $(1,-\frac{1}{2})^\top$.}
    \label{fig:Limit_Points}
\end{figure}
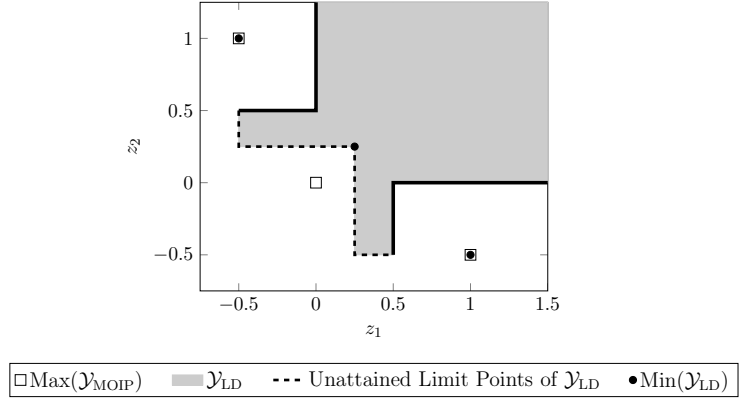

\begin{remark}
The set $\mathcal{Y}_{\mathrm{LD}}$ may be non-convex, disconnected, and neither open nor closed; this is illustrated in Example \ref{ex:P_not_open} and Figure \ref{fig:Limit_Points}.
\end{remark}

\begin{example}\label{ex:P_not_open}
Consider the MOIP
\begin{equation*} \begin{aligned}
\max \enskip & \begin{bmatrix} 1 & -\frac{1}{2} \\ -\frac{1}{2} & 1\end{bmatrix}\begin{bmatrix} x_1\\ x_2\end{bmatrix} \quad 
\text{s.t.} \enskip
x_1 + x_2 \leqq 1, \ x_1,x_2 \in \{0,1\}.
\end{aligned} \end{equation*}

\noindent
As in Example \eqref{ex:LR_better_everywhere}, we dualize the linear constraints. For a fixed 
$\Lambda$, 
this yields the following Lagrangian relaxation.
\begin{equation}
\label{eq:LR_geomP}
\begin{aligned}
\max \enskip & 
\begin{bmatrix} 1 & -\frac{1}{2} \\ -\frac{1}{2} & 1\end{bmatrix} 
\begin{bmatrix} x_1\\ x_2\end{bmatrix}
+ \begin{bmatrix} \lambda_1\\ \lambda_2 \end{bmatrix} \left(1 - \begin{bmatrix} 1 & 1 \end{bmatrix} \begin{bmatrix} x_1\\ x_2\end{bmatrix} \right) \quad
\text{s.t.} \enskip x_1, x_2 \in \{0,1\}.
\end{aligned} 
\end{equation}

For each $\Lambda$, let $\mathcal{Y}_{\mathrm{LR}(\Lambda)}$ be the set of feasible objective values of \eqref{eq:LR_geomP}. Enumerating the feasible objective values, we have
\[
\mathcal{Y}_{\mathrm{LR}(\Lambda)} = \left\{ (\lambda_1, \lambda_2)^\top, \Big(1,-\frac{1}{2}\Big)^\top,\Big(-\frac{1}{2},1\Big)^\top, \Big(\frac{1}{2}-\lambda_1, \frac{1}{2}-\lambda_2\Big)^\top \right\}.
\]
Then, the nondominated points of \eqref{eq:LR_geomP} are given by $\Max(\mathcal{Y}_{\mathrm{LR}(\Lambda)})$.
Considering various sub-cases, we derive the following description of $\Max(\mathcal{Y}_{\mathrm{LR}(\Lambda)})$.
{
\medmuskip=0mu
\thinmuskip=0mu
\nulldelimiterspace=0pt
\[
\Max(\mathcal{Y}_{\mathrm{LR}(\Lambda)})
\hspace{-2pt}=\hspace{-2pt} 
\begin{cases}
\{ (\lambda_1, \lambda_2) \}, & \lambda_1, \lambda_2 \geq 1, \\
\{ (\lambda_1, \lambda_2), (-\frac{1}{2},1) \}, & \lambda_1 \geq 1, 0 \leq \lambda_2 < 1, \\
\{ (\lambda_1, \lambda_2), (1,-\frac{1}{2}) \}, & 0 \leq \lambda_1< 1, \lambda_2 \geq 1, \\
\{ (\lambda_1, \lambda_2), (-\frac{1}{2},1), (1,-\frac{1}{2}) \}, & \frac{1}{4} \leq \lambda_1, \lambda_2 < 1, \\
\{ (-\frac{1}{2},1), (1,-\frac{1}{2}), (\frac{1}{2}-\lambda_1, \frac{1}{2}-\lambda_2) \}, & 0 \leq \lambda_1, \lambda_2 < \frac{1}{4}, \\
\{ (\lambda_1, \lambda_2), (-\frac{1}{2},1), (1,-\frac{1}{2}), (\frac{1}{2}-\lambda_1, \frac{1}{2}-\lambda_2) \}, \hspace{-6pt} & 0 \leq \lambda_1 < \frac{1}{4}, \frac{1}{4} \leq \lambda_2 < 1, \\
\{ (\lambda_1, \lambda_2), (-\frac{1}{2},1), (1,-\frac{1}{2}), (\frac{1}{2}-\lambda_1, \frac{1}{2}-\lambda_2) \}, \hspace{-6pt} & \frac{1}{4} \leq \lambda_1 < 1, 0 \leq \lambda_2 < \frac{1}{4}.
\end{cases}
\]
}
\noindent The set $\mathcal{Y}_{\mathrm{LD}}$ is plotted in Figure \ref{fig:Limit_Points}.

We first show that $\mathcal{Y}_{\mathrm{LD}}$ is not open. 
For $\Lambda = (0,0)^\top$, we have $(1, -\frac{1}{2})^\top \in\Max(\mathcal{Y}_{\mathrm{LR}(\Lambda)}) \subseteq \mathcal{Y}_{\mathrm{LD}}$. Consider an arbitrary $\epsilon > 0$ and suppose $(1 - \epsilon,-\frac{1}{2}-\epsilon)^\top \in \mathcal{Y}_{\mathrm{LD}}$.
Then, there exists $\Lambda \geqq 0$ such that 
$(1 - \epsilon,-\frac{1}{2}-\epsilon)^\top$ is a nondominated element of the set
$\mc{Y}_{\mathrm{LR}(\Lambda)}$. 
However, $(1 - \epsilon, -\frac{1}{2}-\epsilon)^\top \leq (1,-\frac{1}{2})^\top$ for any such $\Lambda$, which is a contradiction. Thus, for each $\epsilon >0$, the ball of radius $\epsilon \sqrt{2}$ centered at $(1,-\frac{1}{2})^\top$ is not contained in $\mathcal{Y}_{\mathrm{LD}}$, so $\mathcal{Y}_{\mathrm{LD}}$ is not open.

Next we show that $\mathcal{Y}_{\mathrm{LD}}$ is not closed either. For this, we show that $(\frac{1}{4},0)^\top$ is a limit point of $\mathcal{Y}_{\mathrm{LD}}$ that is not contained in $\mathcal{Y}_{\mathrm{LD}}$. For all $0 < \epsilon <\frac{1}{4}$, setting $\lambda_1 = \frac{1}{4} + \epsilon$ and $\lambda_2 = 0$ gives
\begin{equation*}
\begin{aligned}
\Max(\mathcal{Y}_{\mathrm{LR}((1/4 +\epsilon,0)^\top)})
&=\left\{\Big(\frac{1}{4} + \epsilon,0\Big)^\top,\Big(1,-\frac{1}{2}\Big)^\top,\Big(-\frac{1}{2},1\Big)^\top, \Big(\frac{1}{4}-\epsilon,\frac{1}{2}\Big)^\top\right\}.
\end{aligned} \end{equation*}

Thus, for all $\epsilon$ small enough, $(\frac{1}{4} + \epsilon, 0)^\top$ is an element of $\mathcal{Y}_{\mathrm{LD}}$. It follows that $(\frac{1}{4},0)^\top$ is a limit point of $\mathcal{Y}_{\mathrm{LD}}$. However, for $(\frac{1}{4}, 0)^\top$ to be a feasible objective value to $LR(\Lambda)$, we must have $\Lambda = (\frac{1}{4}, 0)^\top$ or $\Lambda = (\frac{1}{4}, \frac{1}{2})^\top$. In either case,
\begin{equation*} \begin{aligned}
\Max(\mathcal{Y}_{\mathrm{LR}(\Lambda)})
&= \left\{\Big(1,-\frac{1}{2}\Big)^\top,\Big(-\frac{1}{2},1\Big)^\top,\Big(\frac{1}{4},\frac{1}{2}\Big)^\top
\right\}.
\end{aligned} \end{equation*}

Hence, $\mathcal{Y}_{\mathrm{LD}}$ has an unattained limit point, and $\mathcal{Y}_{\mathrm{LD}}$ is not closed. 
\qed \end{example}

Because $\mathcal{Y}_{\mathrm{LD}}$ is not guaranteed to be closed, one possible avenue is to define the Lagrangian dual via the closure of $\mathcal{Y}_{\mathrm{LD}}$, denoted by $\mathrm{cl}(\mathcal{Y}_{\mathrm{LD}})$. However, this may not be fruitful because points in $\Min(\mathrm{cl}(\mathcal{Y}_{\mathrm{LD}}))$ do not necessarily provide upper bounds on $\Max(\mathcal{Y}_{\mathrm{MOIP}})$. 
To see this in Example \ref{ex:P_not_open}, consider the point $(\frac{1}{4},-\frac{1}{2})^\top$ that lies on the boundary of $ \mathrm{cl}(\mathcal{Y}_{\mathrm{LD}})$ but is not contained in $\mathcal{Y}_{\mathrm{LD}}$, and is dominated by the primal nondominated point $(1,-\frac{1}{2})^\top$.

This breakdown at unattained limit points of $\mathcal{Y}_{\mathrm{LD}}$ can occur even for more structured problems. For example, consider the primal problem
\begin{equation*} \label{eq:primal_ideal}
    \max\left\{\begin{bmatrix}x_1\\x_2 \end{bmatrix}\ \Big\vert \  x_1 + x_2 \leq 2, x_1,x_2 \in \{0,1\} \right\}.
\end{equation*}
\noindent It has a unique efficient solution at $(1,1)^\top$, but dualizing the constraints $x_1 \leq 1$ and $x_2 \leq 1$ and taking $\Lambda_\epsilon = \begin{bmatrix} 2 & 2\\ 0 & 1-\epsilon \end{bmatrix}$ yields 
$y_{\epsilon} = \begin{bmatrix} 0\\ 1+ \epsilon \end{bmatrix} \in \Max(\mathcal{Y}_{\mathrm{LR}(\Lambda_\epsilon)})$ for small $\epsilon > 0$. However, the limit point $(0,1)^\top$ is dominated by $(1,1)^\top$.

Owing to the geometry of $\mathcal{Y}_{\mathrm{LD}}$, we note that points in $\mathcal{Y}_{\mathrm{LD}}$ may provide better bounds than $\Min(\mathcal{Y}_{\mathrm{LD}})$.
For instance, in Example \ref{ex:P_not_open}, $(\frac{1}{4}, \frac{1}{4})^\top$ is the only upper bound in $\Min(\mathcal{Y}_{\mathrm{LD}})$ on the nondominated point $(0,0)^\top$. However, points $(\frac{1}{4} + \epsilon,0)^\top$ and $(0,\frac{1}{4} + \epsilon)^\top$ obtained via Lagrangian relaxation provide tighter upper bounds (for small $\epsilon$). These observations suggest that considering (a finite set of) Lagrangian relaxations may yield a better upper bound set than one obtained by restricting our attention to $\Min(\mathcal{Y}_{\mathrm{LD}})$. A numerical illustration of the quality of (an approximation to) $\Min(\mathcal{Y}_{\mathrm{LD}})$ as an upper bound set is given in Section \ref{sec:numerical}. Note that any approximation to $\mathcal{Y}_{\mathrm{LD}}$
that considers finitely many values of $\Lambda$ will result in a closed set that is bounded below, so that its nondominated set is guaranteed to be nonempty.


While it may be difficult to obtain an explicit description of $\mathcal{Y}_{\mathrm{LD}}$, we can approximate it through its relationship with the set $\mathcal{Y}_{\mathrm{MOIP}}$.
Proposition \ref{prop:P_sub_int_union} shows that for any $y \in \mathcal{Y}_{\mathrm{MOIP}}$, $\mathrm{cl}(\mathcal{Y}_{\mathrm{LD}})$ is contained in the union of the half-spaces $\{z \in \mathbb{R}^{k}\ \vert \ z_i \geq y_i\}$, $i=1, \ldots, k$. 


\begin{proposition}\label{prop:P_sub_int_union}
If $z \in \mathrm{cl}(\mathcal{Y}_{\mathrm{LD}})$, then for each $y \in \mathcal{Y}_{\mathrm{MOIP}}$, $z_i \geq y_i$ for some $i$. 
That is, $\mathrm{cl}(\mathcal{Y}_{\mathrm{LD}}) \subseteq \bigcap\limits_{y \in \mathcal{Y}_{\mathrm{MOIP}}}\bigcup\limits_{i = 1}^{k}\{z \in \mathbb{R}^{k}\ \vert \ z_i \geq y_i\}.$ 
\end{proposition} 

\proof 
Consider a fixed $y \in \mathcal{Y}_{\mathrm{MOIP}}$ and $z \in \mathcal{Y}_{\mathrm{LD}}$. By Proposition \ref{prop_weak_dual}, there exists $i \in \{1, \ldots, k\}$ for which $z_i \geq y_i$. That is, $\mathcal{Y}_{\mathrm{LD}} \subseteq \{z \in \mathbb{R}^{k}\ \vert \ z_i \geq y_i\}$ for some $i=1, \ldots, k$. As $y \in \mathcal{Y}_{\mathrm{MOIP}}$ was arbitrary, we have $\mathcal{Y}_{\mathrm{LD}} \subseteq \cap_{y \in \mathcal{Y}_{\mathrm{MOIP}}} \cup_{i=1}^k \{z \in \mathbb{R}^{k}\ \vert \ z_i \geq y_i\}$.

To see that the containment extends to $\mathrm{cl}(\mathcal{Y}_{\mathrm{LD}})$, note that for each $y \in \mathcal{Y}$, $\bigcup_{i = 1}^{k}\{z \in \mathbb{R}^{k}\ \vert \ z_i \geq y_i\}$ is a finite union of closed sets and is therefore closed. Then, $\bigcap_{y \in \mathcal{Y}_{\mathrm{MOIP}}}\bigcup_{i = 1}^{k}\{z \in \mathbb{R}^{k}\ \vert \ z_i \geq y_i\}$ is also closed so that $\mathcal{Y}_{\mathrm{LD}} \subseteq \mathrm{cl}(\mathcal{Y}_{\mathrm{LD}}) \subseteq \bigcap_{y \in \mathcal{Y}_{\mathrm{MOIP}}}\bigcup_{i = 1}^{k}\{z \in \mathbb{R}^{k}\ \vert \ z_i \geq y_i\}$.
\qed
\endproof

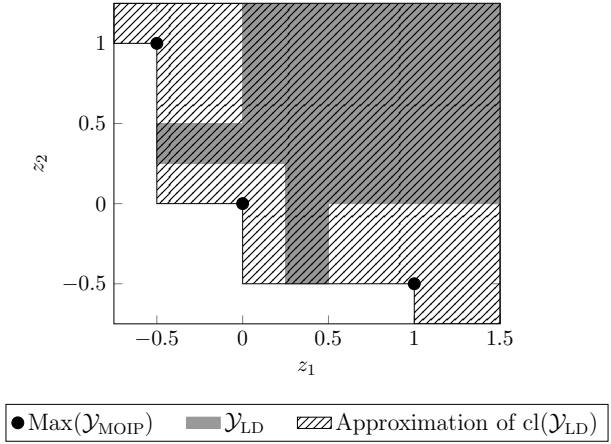
\begin{figure}
    \centering
    \resizebox{0.5\textwidth}{!}{
\begin{tikzpicture}[scale = 0.8]
\begin{axis}[
xmin = {-0.75},
xmax = {1.5},
ymin = {-0.75},
ymax = {1.25},
legend style={at={(0.5,-0.25)},anchor=north},
legend style={/tikz/every even column/.append style={column sep=0.5cm}},
legend columns = -1,
xlabel = {$z_1$},
ylabel = {$z_2$}
]

\addplot[only marks, mark size = 3] coordinates {(0,0) (-0.5,1) (1, -0.5)};
\addplot[fill, black!40!white, area legend] coordinates{
(-0.5,0.25)
(-0.5,0.5)
(0,0.5)
(0,1.25)
(1.5,1.25)
(1.5,0)
(0.5,0)
(0.5,-0.5)
(0.25,-0.5)
(0.25,0.25)
} 
--cycle;



\addplot[fill, pattern = north east lines, area legend] coordinates{
(-0.75,1.25)
(-0.75,1)
(-0.5,1)
(-0.5,0)
(0,0)
(0,-0.5)
(1,-0.5)
(1,-0.75)
(1.5,-0.75)
(1.5,1.25)
}
--cycle;





\legend{{\large 
$\Max(\mathcal{Y}_{\mathrm{MOIP}})$},
{\large $\mathcal{Y}_{\mathrm{LD}}$},
{\large Approximation of $\mathrm{cl}(\mathcal{Y}_{\mathrm{LD}})$}
}

\end{axis}
\end{tikzpicture}
}
    \caption{The set $\mathcal{Y}_{\mathrm{LD}}$ for the MOIP in Example \ref{ex:P_not_open} as well as the approximation of its closure given by Proposition \ref{prop:P_sub_int_union}. Elements of $\mathrm{cl}(\mathcal{Y}_{\mathrm{LD}})$ are weakly dominated by elements of $\Max(\mathcal{Y}_\mathrm{{MOIP}})$. However, for every  $y \in \mathcal{Y}_{\mathrm{MOIP}}$ and $z \in \mathrm{cl}(\mathcal{Y}_{\mathrm{LD}})$, there is an objective $i$ for which $y_i \leqq z_i$.}
    \label{fig:P_approx}
\end{figure}

Figure \ref{fig:P_approx} illustrates the superset described in Proposition \ref{prop:P_sub_int_union} on Example \ref{ex:P_not_open}. Some limit points of $\mathcal{Y}_{\mathrm{LD}}$ are dominated by primal nondominated points. Nonetheless, every element of $\mathrm{cl}(\mathcal{Y}_{\mathrm{LD}})$ has at least one component that exceeds the corresponding objective value of a primal feasible solution. 

\subsection{Properties of the Lagrangian dual}
\label{sec:LagDualProps}

We now proceed to establish properties of the Lagrangian dual \eqref{eq:Lagrange_dual}.

\begin{corollary}[Weak Duality for \eqref{eq:Lagrange_dual}]\label{cor:weak_LD}
If $x$ is feasible to \eqref{eq:MOIP} and $\Lambda \in \reals^{k \times m_1}$, then $\{ Cx \} \preceqq  \Max(\mathcal{Y}_{\mathrm{LR}(\Lambda)})$.
\end{corollary}
\proof
By Proposition \ref{prop_weak_dual}, we have $\{Cx\} \preceqq \Max(\mathcal{Y}_{\mathrm{MOIP}}) \preceqq \Max(\mathcal{Y}_{\mathrm{LR}(\Lambda)})$.
\qed \endproof

Note that Corollary \ref{cor:weak_LD} does not imply $\Max(\mathcal{Y}_{\mathrm{MOIP}}) \preceqq \Min(\mathcal{Y}_{\mathrm{LD}})$. As discussed in Remark \ref{rem:set_order_fail_min}, for $S,T \subseteq \reals^k$ with $S \preceqq T$, we do not necessarily have $S \preceqq \Min(T)$. 
We can only assert that elements of $\Max(\mathcal{Y}_{\mathrm{MOIP}})$ do not dominate elements of $\Min(\mathcal{Y}_{\mathrm{LD}})$.   

Theorem \ref{thm:LD_as_LP} states that the optimal value of the Lagrangian dual of an IP can be obtained by solving an LP. An analogous result does not hold for the multiobjective problem.
That is, the Lagrangian dual of an MOIP cannot be posed as an MOLP in general. This is illustrated in Example \ref{ex:LR_better_everywhere}. 
The set $\Min(\mc{Y}_{\mathrm{LD}})$ consists of three isolated points which can never be obtained from a single MOLP \edit{because the nondominated set of an MOLP is connected}. Nonetheless, Theorem \ref{thm:MOLP_LD_Bound} uses an MOLP to derive an upper bound on the dual nondominated points. Recall that if $A^1$ and $A^2$ are the sub-matrices of $A$ corresponding to the complicating and simple constraints respectively, then $Q = \{x \in \mathbb{Z}^n_\geqq \ \vert \ A^2 x \leqq b^2\}$.

\begin{theorem} \label{thm:MOLP_LD_Bound} 
Let $\mathcal{Y}_{\mathrm{LDLP}}$ be the set of feasible objective values of the MOLP
\begin{equation}
\label{eq:MOLP-LD}
\max\ \{Cx\ \vert \ A^1x \leqq b^1, \ x \in \text{\emph{conv}}(Q)\}.
\end{equation}
Then, $\Min(\mc{Y}_{\mathrm{LD}}) \preceqq \Max(\mathcal{Y}_{\mathrm{LDLP}})$.
\end{theorem}

\proof 
Because $\text{conv}(Q)$ is a polyhedron, there is a matrix $B$ and a vector $d$ such that $\text{conv}(Q) = \{x \in \reals^n_\geqq\ \vert \ Bx \leqq d\}$. Strong Lagrangian duality for MOLPs (Lemma \ref{lem:MOLP_LD_STRONG}) implies that
\begin{align}
\nonumber \Max(\mathcal{Y}_{\mathrm{LDLP}})
& = \max\ \{Cx \ \vert \ A^1x \leqq b^1, x\ \in \text{conv}(Q)\}  \\
\nonumber &= \min_{\Lambda,\Gamma \geqq 0}\ \max_{x \in \reals^n}\ ( Cx + \Lambda(b^1 - A^1x) + \Gamma(d - Bx) )\\
\label{eqn:LD_MOLP_1} & \suceqq \min_{\Lambda \geqq 0} ~\max_{x \in \text{conv}(Q)}\ (Cx + \Lambda(b^1- A^1x) ) \\
\label{eqn:LD_MOLP_2} &\suceqq \min_{\Lambda \geqq 0} ~\max_{x \in Q}\ (Cx + \Lambda(b^1- A^1x))\\
\nonumber &\suceqq \Min\ (\mc{Y}_{\mathrm{LD}}).
\end{align}

The inequality in \eqref{eqn:LD_MOLP_1} holds because for a fixed $\Lambda \geqq 0$ and for every $\Gamma \geqq 0$, $x \in \text{conv}(Q)$, we have
$Cx + \Lambda(b^1 - A^1x) + \Gamma(d - Bx) \geqq Cx + \Lambda(b^1 - A^1x)$. In particular, 
\[
 \max\limits_{x \in \reals^n}\{Cx + \Lambda(b^1 - A^1x) + \Gamma(d - Bx)\}\suceqq  \max\{Cx + \Lambda(b^1 - A^1x) \ \vert \ Bx \leqq d\}.
\]
Further, \eqref{eqn:LD_MOLP_2} holds because for any $\Lambda \geqq 0$, if $x \in Q$ is an efficient solution for $\max\limits_{x \in Q}\ (Cx + \Lambda(b^1- A^1x))$, then it is either also an efficient solution for $\max\limits_{x \in \text{conv}(Q)}\ (Cx + \Lambda(b^1 - A^1x))$, or an interior point of $\text{conv}(Q)$. In the latter case, $\{Cx +\Lambda(b^1 - A^1)x\} \preceqq \max\limits_{x \in \text{conv}(Q)}\ (Cx + \Lambda(b^1 - A^1x)) $. Thus, for each $\Lambda \geqq 0$, 
\begin{equation*} 
\begin{aligned}
\max_{x \in \text{conv}(Q)}\ (Cx + \Lambda(b^1- A^1x)) \ \suceqq \ \max_{x \in Q} \ (Cx + \Lambda(b^1 - A^1x)). 
\end{aligned} 
\tag*{\qed}
\end{equation*} 
\endproof

Corollary \ref{Cor_LD_better_MOLP} establishes that the Lagrangian dual provides a tighter upper bound than that given by the continuous relaxation \eqref{eq:MOLP}.

\begin{corollary} \label{Cor_LD_better_MOLP} Let $\mathcal{Y}_{\mathrm{MOLP}}$ be the set of feasible objective values for \eqref{eq:MOLP}. Then, $\Min(\mathcal{Y}_{\mathrm{LD}}) \preceqq \Max(\mathcal{Y}_{\mathrm{MOLP}})$.
\end{corollary}

\proof 
The feasible region of \eqref{eq:MOLP} contains the feasible region of \eqref{eq:MOLP-LD}. Therefore,
$
\Max(\mathcal{Y}_{\mathrm{MOLP}}) \suceqq \Max(\mathcal{Y}_{\mathrm{LDLP}}) \suceqq \Min(\mathcal{Y}_{\mathrm{LD}})
$.
\qed \endproof

Corollary \ref{Cor_LD_better_MOLP} implies that if $x^*$ and $\tilde{x}$ are efficient solutions to \eqref{eq:MOIP} and its MOLP relaxation respectively and $y \in \Min(\mathcal{Y}_{\mathrm{LD}})$ is such that $Cx^*,C\tilde{x}$, and $y$ are all comparable, then $Cx^{*} \leqq y \leqq C \tilde{x}$. Example \ref{ex:LR_better_everywhere} illustrates that both the inequalities can be strict. 

Theorem \ref{thm:MOLP_LD_Bound} gives a loose upper bound on the Lagrangian dual. In the remainder of this section, we investigate relationships that hold with equality. To do so, we solve a scalarized problem and apply results from single-objective duality.

\begin{theorem}\label{thm_scalarized_LDLP}
For all $\mu \in \mathbb{R}^k_{>}$, 
\begin{equation*} 
\begin{aligned}
\min_{\Lambda \geqq 0}\ \max_{x \in Q}\ \mu^\top Cx + \mu^\top \Lambda(b^1-A^1x) = \max_{x \in \text{\emph{conv}}(Q)} \{ \mu^\top Cx\ \vert \ A^1x \leqq b^1 \}.
\end{aligned} 
\end{equation*}
\end{theorem}

\proof 
Consider the single-objective IP 
\begin{equation}\label{eq_IP_scalarized}
\begin{aligned}
 \max \ & \mu^\top C x \quad
 \text{s.t.}\enskip A^1x \leqq b^1, \ x \in Q.
 \end{aligned}
\end{equation}
The Lagrangian dual of \eqref{eq_IP_scalarized} is
\begin{equation}\label{eq_LD1_scalarized_LDLP}
\min_{u \in \mathbb{R}^m_\geqq}\ \max_{x \in Q}\ \mu^\top Cx + u^\top (b^1 - A^1x).
\end{equation}
For any $u \in \mathbb{R}^{m_1}_\geqq$, we can write $u^\top = \mu^\top \Lambda$ for some $\Lambda \in \mathbb{R}^{k\times m_1}_{\geqq}$. Conversely, if $\mu\in \reals^k_{>}$ and $\Lambda \in \reals^{k \times m_1}_{\geqq}$, then $(\mu^\top \Lambda)^\top \in \reals^{m_1}_\geqq$. Thus, we can rewrite \eqref{eq_LD1_scalarized_LDLP} as 
\begin{equation}\label{eq_LD2_scalarized_LDLP}
\min_{\Lambda \in \reals^{k \times m_1}_\geqq}\ \max_{x \in Q}\ \mu^\top Cx + \mu^\top \Lambda(b^1 - A^1x).
\end{equation}

Moreover, Theorem \ref{thm:LD_as_LP} 
implies that the Lagrangian dual \eqref{eq_LD1_scalarized_LDLP} to \eqref{eq_IP_scalarized} has the same optimal value as the LP

\begin{equation}\label{eq_LP_scalarized}
\begin{aligned}
 \max \enskip & \mu^\top C x \quad \text{s.t.}\enskip  A^1x \leqq b^1, \ x \in \text{conv}(Q).
 \end{aligned}
\end{equation}
Therefore,
\begin{equation*} 
\begin{aligned}
\min_{\Lambda \geqq 0}~\max_{x \in Q}~ \mu^\top Cx + \mu^\top \Lambda(b^1-A^1x) = \max_{x \in \text{conv}(Q)}~ \left\{\mu^\top Cx\ \vert \ A^1x \leqq b^1 \right\}.
\end{aligned}
\tag*{\qed}
\end{equation*}
\endproof

Thus, single-objective LPs can solve the scalarized dual problem.

\subsection{Strong Lagrangian Duality}

A dual problem to \eqref{eq:MOIP} is strong at a solution $x$ if there exists a point $y \in \Min(\mathcal{Y}_{\mathrm{LD}})$ such that $y = Cx$.
In this subsection, we seek conditions on \eqref{eq:MOIP} under which the Lagrangian dual is strong. Theorem \ref{thm:MOLP_LD_Bound} used an MOLP to prescribe an upper bound for $\Min(\mathcal{Y}_{\mathrm{LD}})$. Theorem \ref{Strong_Dual_sufficient} uses this MOLP to derive a sufficient condition for strong duality.

\begin{theorem}\label{Strong_Dual_sufficient}
Let $x^*$ be an efficient solution of \eqref{eq:MOIP} such that $Cx^* \leqq y$ for some $y \in \Min(\mathcal{Y}_{\mathrm{LD}})$. 
The Lagrangian dual \eqref{eq:Lagrange_dual} is strong at $x^*$ 
if there exists a vector $\alpha \in \mathbb{R}^k_{>}$ such that 
\begin{equation}
\label{eq_StrongDual_suff}
\alpha^\top C(x^*-x) \leqq 0 \text{ \emph{for all} } x \in \text{\emph{conv}}(Q) \cap \{x \in \mathbb{R}^n\ \vert \ A^1x \leqq b^1\}.
\end{equation} 

\end{theorem}

\proof 
Condition \eqref{eq_StrongDual_suff} holds if and only if $x^*$ is an efficient solution to the MOLP
$\max \{Cx\ \vert \ A^1 \leqq b^1, x \in \text{conv}(Q)\}$ \cite[Theorem 4.2.6(i)]{luc_2016_Book}.
Moreover, Theorem \ref{thm:MOLP_LD_Bound} implies that $\Min(\mathcal{Y}_{\mathrm{LD}}) \preceqq \max \{Cx\ \vert \ A^1 \leqq b^1, x \in \text{conv}(Q)\}$ so that if $y \in \Min(\mathcal{Y}_{\mathrm{LD}})$ is comparable to $Cx^*$, then $y \leqq Cx^*$. Corollary \ref{cor:weak_LD} then implies that $y \not \leq Cx^*$. Thus, $y = Cx^*$.
\qed \endproof

Restricting our attention to supported solutions, we next derive conditions for strong Lagrangian duality that are independent of the objective function. These results are analogous to the single-objective case (see Theorem \ref{thm:Single_objective_LD_iff}).

\begin{theorem}\label{thm:MO_Lag_iff_1}
If the Lagrangian dual is strong at supported efficient solutions for all matrices $C$, then 
\begin{equation}\label{eq:FR_LAG}\
\text{\emph{conv}}(Q \cap \{x \in \mathbb{R}^n\ \vert \ A^1x \leqq b^1\}) = \text{\emph{conv}}(Q) \cap \{x \in \mathbb{R}^n\ \vert \ A^1x \leqq b^1\}.
\end{equation}
\end{theorem}

\proof 
\edittwo{ 
By Theorem \ref{thm:Single_objective_LD_iff}, it suffices to show that under the hypothesis of Theorem \ref{thm:MO_Lag_iff_1}, strong Lagrangian duality holds for any single-objective IP with the same feasible region as \eqref{eq:MOIP}, i.e, $Q \cap \{x \in \mathbb{R}^n\ \vert \ A^1x \leqq b^1\}$. 

Let $c \in \reals^n$ be an arbitrary cost vector and consider the cost matrix $C = \begin{bmatrix} c & 0 & \dots & 0\end{bmatrix}^\top \in \reals^{k \times n}$. Set $\mu \in \reals^k_{>}$ to be the vector of all ones so that $\mu^\top C = c^\top$. Let $x^*$ be a supported efficient solution to \eqref{eq:MOIP} with supporting vector $\mu$. Then $x^*$ is an efficient solution to the single-objective IP $\max \{c^\top x \ \vert \ A^1x \leqq b^1, x \in Q\}$.  

By the hypothesis, there exists $y \in \Min(\mathcal{Y}_{\mathrm{LD}})$ such that $Cx^* = y$. So, $y = C\hat{x} + \hat{\Lambda}(b^1 - A^1\hat{x})$ for some $\hat{x} \in Q$ and $\hat{\Lambda} \geqq 0$. Then,
\[
c^\top x^* = \mu^\top Cx^* = \mu^\top y = \mu^\top C\hat{x} + \mu^\top \hat{\Lambda}(b^1 - A^1\hat{x}) = c^\top \hat{x} + \hat{\lambda}(b^1 - A^1\hat{x}),\]
where $\hat{\lambda} = \mu^\top \hat{\Lambda} \in \reals^{m_1}_{\geqq}$. Thus, the Lagrangian dual of the IP $\max \{c^\top x \ \vert \ A^1x \leqq b^1, x \in Q\}$ is strong. Because $c$ was arbitrary, Theorem \ref{thm:Single_objective_LD_iff} implies that \eqref{eq:FR_LAG} holds. 
\qed\endproof
}

An example to illustrate Theorem \ref{thm:MO_Lag_iff_1} is provided in Appendix \ref{sec:ex:LD_not_strong_supp}, where we present an MOIP whose feasible region does not satisfy condition \eqref{eq:FR_LAG}, and show that the Lagrangian dual for this problem is not strong at a supported solution.
A modified converse to Theorem \ref{thm:MO_Lag_iff_1} is derived in Theorem \ref{thm:MO_Lag_iff_2}, but we first show that optimal objective values for the scalarized dual problem lift to feasible objective values for the multiobjective dual. 

\begin{lemma}\label{lem:scalarized_dual}
For any $\mu \in \reals^k_{>}$, let
\[
w = \min_{\Lambda \in \reals^{k \times m_1}_\geqq}\ \max_{x \in Q}\ \mu^\top Cx + \mu^\top \Lambda(b^1-A^1x),\]
\noindent be the optimal objective value for the scalarized dual problem. Then, there exists $y \in \mathcal{Y}_{\mathrm{LD}}$ such that $\mu^\top y = w$.
\end{lemma}

\proof 
Fix $\mu \in \reals_{>}^k$.
Let $w = \min\limits_{\Lambda \in \reals^{k \times m_1}_\geqq} \max\limits_{x \in Q}\ \mu^\top Cx + \mu^\top \Lambda(b^1-A^1x)$. Then, there exists $\Lambda^* \in \reals^{k\times m_1}_\geqq$ and $x^* \in Q$ such that $x^*$ is optimal for the inner maximization and 
$w = \mu^\top y^*$ where $y^* = Cx^* + \Lambda^*(b^1-A^1x^*)$.
Because $x^*$ is optimal to the scalarized problem, it is a supported efficient solution for the relaxation $\mathrm{LR(\Lambda^*)}$. Hence, $y^* \in \Max(\mc{Y}_{\mathrm{LR}(\Lambda^*)}) \subseteq \mathcal{Y}_{\mathrm{LD}}$.
\qed \endproof

\begin{theorem}\label{thm:MO_Lag_iff_2}
Suppose \eqref{eq:FR_LAG} holds. Given an objective matrix $C$ and a supporting vector $\mu \in \reals^k_{>}$, let $x^*$ be a corresponding supported solution to \eqref{eq:MOIP}. If there exists $y \in \Min(\mathcal{Y}_{\mathrm{LD}})$ that is comparable to $Cx^*$, then the Lagrangian dual is strong at $x^*$. 
\end{theorem}

\proof
Let $\mu \in \reals^k_{>}$ and $C \in \reals^{k \times m}$ be arbitrary. If \eqref{eq:FR_LAG} holds, then we have by Theorem \ref{thm:Single_objective_LD_iff} that 
\begin{equation*} \begin{aligned}\max_{x \in Q}\{\mu^\top Cx\ \vert \ A^1x\leqq b^1\} = \min_{\Lambda \geqq 0} \max_{x \in Q} \mu^\top Cx + \mu^\top \Lambda(b^1 - A^1x).\end{aligned} \end{equation*}

Let $x^*$ be an optimal solution to the scalarized problem 
\begin{equation*} \begin{aligned}
 \max_{x \in Q}\{\mu^\top Cx\ \vert \ A^1x \leqq b^1\},
\end{aligned} \end{equation*}
\noindent such that there exists $\tilde{y} \in \Min(\mathcal{Y}_{\mathrm{LD}})$ that is comparable with $Cx^*$. Corollary \ref{cor:weak_LD} implies that $Cx^* \leqq \tilde{y}$. Lemma \ref{lem:scalarized_dual} implies that there is a $y \in \mathcal{Y}_{\mathrm{LD}}$ such that $\mu^\top Cx^* = \mu^\top y$. But then, because 
\begin{equation*} 
\begin{aligned}
\mu^\top y = \min_{\Lambda \geqq 0} \max_{x \in Q} \mu^\top Cx + y^\top \Lambda(b^1 - A^1x) = \max_{x \in \text{conv}(Q)}\{\mu^\top Cx\ \vert \ A^1x \leqq b^1\},
\end{aligned} 
\end{equation*}
\noindent 
and $\{\tilde{y}\} \preceqq \max_{x \in \text{conv}(Q)}\{Cx\ \vert \ A^1x \leqq b^1\}$ by Theorem \ref{thm:MOLP_LD_Bound}, it follows that
\begin{equation*} 
\begin{aligned}
\mu^\top \tilde{y} \leqq \max_{x \in \text{conv}(Q)}\{\mu^\top Cx\ \vert \ A^1x \leqq b^1\}.
\end{aligned} 
\end{equation*}
\noindent 
Then, 
\begin{equation*} \begin{aligned}\mu^\top y = \mu^\top Cx^* \leqq \mu^\top \tilde{y} \leqq \mu^\top y.
\end{aligned} 
\end{equation*}
This implies that $\mu^\top Cx^* = \mu^\top \tilde{y}$. So, $\mu^\top (\tilde{y} -Cx^*) = 0$ and $\tilde{y} - Cx^*$ has all nonnegative entries, so that $\tilde{y} -Cx^* = 0$. Thus, strong duality holds at $x^*$.
\qed \endproof

In the single-objective case, every efficient solution is a supported solution and the extended real line is totally ordered. Then, Theorem \ref{thm:MO_Lag_iff_2} is the converse of Theorem \ref{thm:MO_Lag_iff_1}, and together, they coincide with Theorem \ref{thm:Single_objective_LD_iff}. 

\begin{remark}
Theorems \ref{thm:MO_Lag_iff_1} and \ref{thm:MO_Lag_iff_2} do not address the strength of the Lagrangian dual at unsupported solutions. Example \ref{ex:LD_not_strong_not_supp} illustrates that even if the feasible region satisfies condition \eqref{eq:FR_LAG} and the Lagrangian dual is strong at supported solutions, it may not be strong at an unsupported solution. 
\end{remark}

\begin{example}\label{ex:LD_not_strong_not_supp}
We revisit the MOIP from Example \ref{ex:P_not_open}.
For this problem, we have $\Max(\mathcal{Y}_{\mathrm{MOIP}})  \allowbreak = \{(1,-\frac{1}{2})^\top, (0,0)^\top, (-\frac{1}{2},1)^\top\}$. For $\mu = (\mu_1, \mu_2)^\top \in \reals^2_{>}$, the scalarized IP is
\begin{equation*} \begin{aligned}
\max \enskip & \Big(\mu_1 - \frac{1}{2}\mu_2\Big)x_1 + \Big(\mu_2 - \frac{1}{2}\mu_1\Big)x_2\\
\text{s.t.} \enskip &
x_1 + x_2 \leqq 1,\\
& x_1,x_2 \in \{0,1\}.
\end{aligned} \end{equation*}
\noindent 
The optimal solution for this IP is $(1,0)^\top$ if $\mu_1 \geq \mu_2$ and $(0,1)^\top$ for $\mu_2 \geq \mu_1$. Thus, $(1,0)^\top$ and $(0,1)^\top$ are supported efficient solutions to the primal problem. Moreover, Example \ref{ex:P_not_open} shows that their corresponding objective vectors $(1,-\frac{1}{2})^\top$ and $(-\frac{1}{2},1)^\top$ are feasible objective values for the Lagrangian dual problem, so that strong duality holds for the supported efficient solutions. 

On the other hand, $(0,0)^\top$ is not a feasible objective value to the dual problem, so that strong duality does not hold for the unsupported solution $(0,0)^\top$, despite the fact that 
\begin{align*}
\text{conv}\big( \{0,1\}^2 \cap \{x \in \mathbb{R}^n \vert x_1 + x_2 \leqq 1\}\big) 
{=} 
\text{conv}\big(\{0,1\}^2\big) \cap \{x \in \mathbb{R}^n \vert  x_1 + x_2 \leqq 1\}. 
\tag*{\qed}
\end{align*}
\end{example}

Theorems \ref{thm:MO_Lag_iff_1} and \ref{thm:MO_Lag_iff_2} analyze strong duality at the supported solutions of \eqref{eq:MOIP}. 
We now consider the unsupported solutions.
Theorem \ref{thm:not_supp_not_strong} derives a sufficient condition under which \eqref{eq:Lagrange_dual} is {\em not} strong at unsupported solutions.

\begin{theorem}\label{thm:not_supp_not_strong}
Let $x$ be an unsupported efficient solution to \eqref{eq:MOIP}. Suppose there exists $\Lambda^* \in \reals^{k\times m_1}_{\geqq}$ and a supported solution $x^*$ to \eqref{eq:Lagrangian_Relaxation} such that 
\begin{enumerate}[label=\roman*)]
    \item $Cx \leq Cx^* + \Lambda^*(b^1-A^1x^*)$, and
    \item for all $\Lambda \in \reals^{k \times m_1}_{\geqq}$, 
$ \{v \in \Max(\mathcal{Y}_{\mathrm{LR}(\Lambda)})\ \vert\ v \leqq Cx^* + \Lambda^*(b^1 - A^1x^*)\}$
is either empty or consists only of supported objective values.
\end{enumerate}
Then, \eqref{eq:Lagrange_dual} is not strong at $x$.
\end{theorem}

\proof
Note that if $y \in \Min(\mathcal{Y}_{\mathrm{LD}})$ satisfies $y = Cx$, then $y \in \Max(\mathcal{Y}_{\mathrm{LR}(\Lambda)}) \cap \{v \in \reals^k\ \vert \ v \leqq Cx^* + \Lambda^*(b^1 - A^1x^*)\}$ for some $\Lambda \in \reals^{k\times m_1}_\geqq$. By hypothesis, $y$ is a supported nondominated point to \eqref{eq:Lagrangian_Relaxation} so that there exists a vector $\mu \in \reals^k_{>}$ such that 
\begin{equation*}
\mu^\top y = \max_{\xi \in Q}\ \mu^\top C\xi + \mu^\top\Lambda(b^1 - A^1 \xi). \end{equation*}

Then, because $y$ is an optimal objective \edit{value} of the single-objective Lagrangian relaxation of $\max\ \{\mu^\top C\xi\ \vert\ \xi \in \mathcal{X}\}$ we have by Proposition \ref{prop:IP_weak_Lag_Dual}

\begin{equation*}\max_{\xi \in \mathcal{X}}\ \mu^\top C\xi \leq \mu^\top y.\end{equation*}

On the other hand, if $x$ is an unsupported efficient solution to \eqref{eq:MOIP}, then $\mu^\top Cx < \max\limits_{\xi \in \mathcal{X}}\ \mu^\top C\xi$. So, $Cx \not = y$. 
\qed\endproof

\subsection{Numerical Illustration}\label{sec:numerical}

As discussed in Section \ref{sec:MODualBoundSets}, bound sets are a key component of search-based algorithms for solving MOIPs. In this section, we derive an upper bound set that approximates the Lagrangian dual and test its performance on two biobjective MOIPs. We compare this set with the upper bound set derived by Ehrgott and Gandibleux in \cite{EHRGOTT2007BoundSets}\footnote{\cite{EHRGOTT2007BoundSets} presents a minimization problem, so the roles of upper and lower bound sets are reversed and accordingly adapted here.}, which amounts to the convex hull of supported nondominated points. As such, we present computational evidence that a Lagrangian dual-based approach can provide a tighter upper bound than that obtained via the convex hull relaxation. Remark \ref{rem:NoRelation_CH_LR} noted that the bound sets due to Lagrangian and convex hull relaxations are mutually incomparable in general, but this section illustrates that (an approximation of) the Lagrangian dual may present a computational advantage.

We test the quality of the bound sets on 100 randomly generated instances of two classes of biobjective problems: a linear assignment problem, and a knapsack problem, each with one additional randomly generated constraint that is subsequently dualized.  
The linear assignment problem consists of 16 binary variables. 
Problem parameters are taken from a binary linear assignment problem appearing in \cite{ulungu1995two}, included as an example file with the Julia-based MOIP solver vOptSolveGeneric \cite{vOptSolver2021}. 
The knapsack problem has 20 binary variables; coefficients for the objectives and constraints are generated randomly in each trial by sampling uniformly at random over the sets $\{1,\ldots, 15\}$ and $\{1,\ldots, 5\}$ respectively.
For both problems, coefficients for the one additional constraint are randomly generated in each trial by sampling uniformly over $\{1,\ldots, 5\}$. All computations were performed in Julia using the vOptSolve package \cite{vOptSolver2021,vOptSolver2017MCDM,vOptSolver2017IFORS} with the GLPK optimizer \cite{glpk}.

To approximate the Lagrangian dual problem, we consider a finite set of Lagrangian relaxations parameterized by multipliers $\Lambda \in \mathcal{M}$, constructed by dualizing the additional constraint for each problem instance. The values of $\Lambda$ in $\mc{M}$ are na\"ively selected as equally spaced gridpoints in $[0,2.5]^2$. The set $\mc{M}$ consists of $51^2$ and $26^2$ values of $\Lambda$ for the linear assignment and knapsack problems respectively. 
Then, $U = \Min \Big( \bigcup\limits_{\Lambda \in \mathcal{M}} \  \Max\left(\mathcal{Y}_{\mathrm{LR}(\Lambda)}\right) \Big)$
is an upper bound set that approximates the Lagrangian bound set $\Min(\mathcal{Y}_{\mathrm{LD}})$ defined \edit{in \eqref{eq:Lagrange_dual}}. 
To determine the quality of this bound, we used a scaled distance 
\begin{equation*}
    d(L,U) = \frac{1}{\gamma} \ \max_{\ell \in L}\ \min_{u \in U}\ \|u - \ell\|_2,
\end{equation*}
where $L$ is a lower bound set consisting of local nadir points between supported solutions \cite{EHRGOTT2007BoundSets}, and $\gamma$ is the average of $\|y \|_2$  for $y \in L \cup U$. A smaller value of $d$ corresponds to a better bound. The set $L$ is the collection of points of the form $(\min( y_1^1,y_1^2), \min(y_2^1,y_2^2))^T$ where the pair $(y^1,y^2)$ ranges over adjacent supported solutions; it is computed as in \cite{EHRGOTT2007BoundSets}.  

\begin{figure}[t]
    \centering
    
    \subfloat[Linear Assignment Problem]{
    \includegraphics[width=0.5\textwidth]{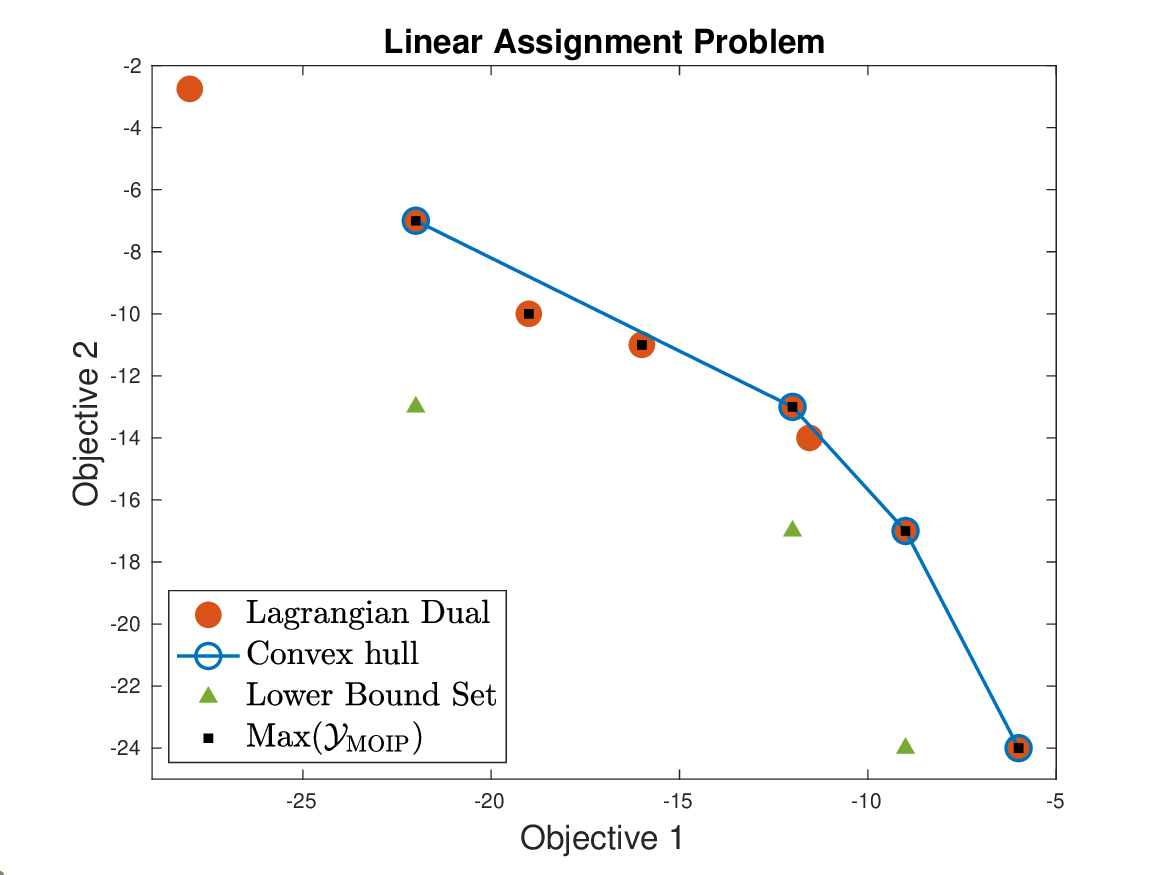}
    }
    \subfloat[Knapsack Problem]{
    \includegraphics[width=0.5\textwidth]{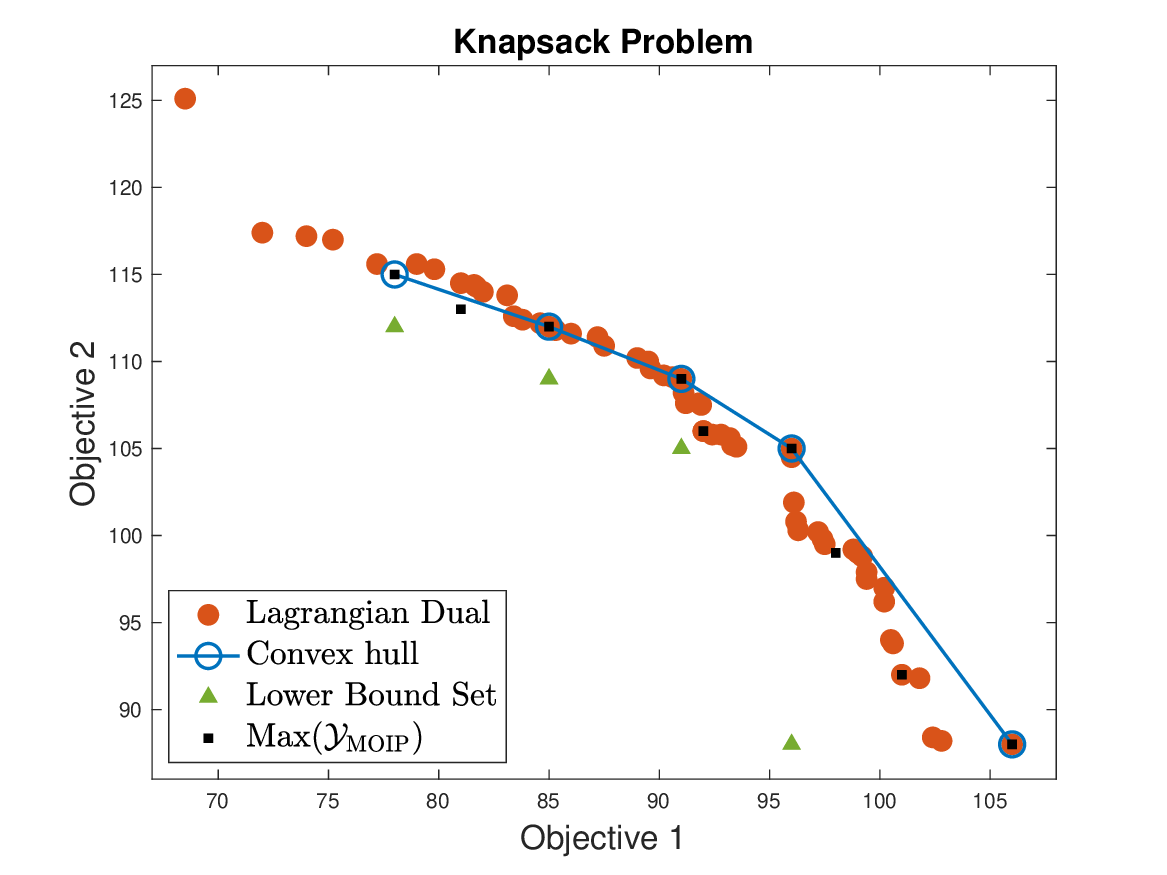}
    }
    \caption{An example of upper bounds computed by an approximation of the Lagrangian dual and the convex hull of the supported nondominated points, as well as a lower bound set derived from the local nadir points of the supported efficient solutions. Note that the Lagrangian upper bound is tight for this instance of the linear assignment problem, but not for the knapsack problem. Due to the presence of unsupported solutions, the convex hull bound is not strong in either instance.}
    \label{fig:ComputationFigures}
\end{figure}

We note that the metric $d$ is a re-scaled version of the measure $\mu_1$ used in \cite{EHRGOTT2007BoundSets}, which is sensitive to outliers in the set $L\cup U$.
Because the set $U$ derived from Lagrangian relaxations can contain points far from the nondominated set of \eqref{eq:MOIP} (as in Figure \ref{fig:ComputationFigures}), this could lead to misleading comparisons and spuriously improve the performance of $U$.
Therefore, we choose the modified scaling parameter $\gamma$ which is less sensitive to far away points. We also report the number of instances for which the upper bound is strong (i.e. $\Max(\mathcal{Y}_{\mathrm{MOIP}}) \subseteq U$). 

We compare the performance of $U$ with $\Max(\mc{Y}_{\mathrm{CH}})$ as computed in \cite[Algorithm 1]{EHRGOTT2007BoundSets}, which recursively solves a sequence of scalarized problems to obtain the supported nondominated points and their convex hull. 
Each single-objective problem is solved exactly---and not approximately as proposed in \cite{EHRGOTT2007BoundSets}---to derive the exact convex hull relaxation.
Again, we measure the quality of this upper bound set using the measure $d$ and the lower bound set $L$ of local nadir points, 
and report the number of instances for which the bound is tight.

\begin{table}[t]
    \centering

    \begin{tabular}{c||c|c|c||c|c|c||}
                & \multicolumn{3}{|c||}{Lagrangian Dual} & \multicolumn{3}{|c||}{Convex Hull}\\ \hline   
        Problem & Mean $d$ & SD $d$ & \# Strong & Mean $d$ & SD $d$ & \# Strong\\
        \hline
        Linear Assignment &3.930& 0.587& 59/100&4.801 & 0.956 & 9/100\\
        Knapsack  & 5.869 & 1.980& 3/100 & 6.384& 1.803 & 9/100\\
    \end{tabular}
    \caption{Mean and standard deviation (SD) of the metric $d$, as well as number of problem instances (out of 100) where the upper bound sets are strong (\#Strong), for bound sets computed by an approximation of the Lagrangian dual and the convex hull of the supported nondominated points.}    \label{tab:LagrangeComputations}
\end{table}
The results are summarized in Table \ref{tab:LagrangeComputations}. 
We find that for both problems, the Lagrangian dual approximation outperforms the convex hull bound with respect to the metric $d$ averaged across the 100 trials. However, the extent of improvement is not uniform across the two problem classes. Specifically, for the linear assignment problem, the standard deviation in $d$ is notably smaller and the Lagrangian dual provides tight bounds in a much larger fraction of trials. In contrast, for the knapsack problem, the convex hull upper bound yields a lower standard deviation in $d$ and provides tight bounds more frequently than the Lagrangian dual. Examples of the bounds for each problem are illustrated in Figure \ref{fig:ComputationFigures}. 

The above numerical experiments are a proof of concept rather than a detailed computational study. The results are promising and point towards several avenues for future investigation. For instance, our approach na\"ively enumerates over a range of Lagrange multiplier matrices $\Lambda$, which limits our ability to scale to larger problems. Future research could explore more sophisticated techniques for searching over the space of multipliers. 
Additionally, in our experiment we use a lower bound set derived from information from the convex hull relaxation of our problem, and the question of similarly using information from the Lagrange relaxation/dual to derive a lower bound set, remains open.  
Although the scope of our experiments is limited, the results are encouraging and highlight the potential of the Lagrangian approach in deriving bound sets for MOIP solution methods.

\section{Superadditive Duality for Multiobjective Integer Programs} \label{sec:superadditive}
In this section, we develop a multiobjective counterpart of the superadditive dual of an IP. 
For $\beta \in \reals^m$, define $\mathcal{X}(\beta) = \{x \in \mathbb{Z}^n_\geqq \ \vert \ Ax \leqq \beta\}$ as the feasible region of \eqref{eq:MOIP} parameterized by the right-hand-side $\beta$. 
We define the value function of an MOIP as
\begin{equation}
\label{eq:Value_Function}
\begin{aligned}
Z(\beta) = \max \enskip Cx \quad
\text{s.t.} \enskip x \in \mathcal{X}(\beta).
\end{aligned}
\end{equation}

Unlike the single-objective case, where the value function maps onto the extended real line, the multiobjective value function $Z$ is {\em multi-valued} in general, and maps to the nondominated set of the MOIP. Moreover, the cardinality of the image set $Z(\beta)$ is not known {\em a priori}. Therefore, the first step towards developing a superadditive dual is to extend the definitions of monotonicity and superadditivity to set-valued functions. Recall from \eqref{eq:setE} that $\mathcal{E}$ is the collection of nonempty subsets of $\reals^k$ whose elements are mutually incomparable, along with $\pm M_{\infty}$.

\begin{definition}[Nondecreasing Function]
A function $F:\reals^m \rightarrow \mathcal{E}$ is nondecreasing (with respect to $\preceqq$) if for all $\beta_1 \leqq \beta_2$, $F(\beta_1) \preceqq F(\beta_2)$.
\end{definition}

\begin{definition}[Superadditive Function]
A function $F: \reals^m \rightarrow \mathcal{E}$ is superadditive (with respect to $\preceqq$) if for all $\beta_1,\beta_2$, $F(\beta_1) + F(\beta_2) \preceqq F(\beta_1+\beta_2)$.
\end{definition}

Here, $F(\beta_1) + F(\beta_2)$ denotes the Minkowski sum of the two sets, defined as $F(\beta_1) + F(\beta_2) = \{ z_1 + z_2\ \vert \ z_1 \in F(\beta_1), z_2 \in F(\beta_2)\}$. 
If $F:\reals^m \rightarrow \mathcal{E}$ is superadditive and nondecreasing with respect to $\preceqq$, then 
$\sum_{j = 1}^{\ell}F(\beta_j) \preceqq F\Big(\sum_{j = 1}^{\ell}\beta_j\Big)$ for any finite sum (by induction on $\ell$). It follows that for any positive integer $\kappa$, we must have $\kappa \cdot F(\beta) \preceqq F(\kappa \cdot \beta)$, where $\kappa \cdot F(\beta)$ is the Minkowski sum of $\kappa$ copies of $F(\beta)$. 

\begin{proposition}
The value function $Z$ is nondecreasing and superadditive with respect to $\preceqq$.
\end{proposition}

\proof 
We first show that $Z$ is nondecreasing. Let $\beta_1 \leqq \beta_2$. Then $\mathcal{X}(\beta_1) \subseteq \mathcal{X}(\beta_2)$. Therefore if $Cx_1 = z_1 \in Z(\beta_1)$, then there is a $z_2 \in Z(\beta_2)$ such that $z_1 \leqq z_2$. Moreover, there is no $z_3 \in Z(\beta_2)$ with $z_3 \leq z_1$, because that would imply $z_3 \leq z_1 \leqq z_2$ which contradicts the nondominance of $z_3$ for the MOIP with RHS $\beta_2$.

Next, we prove that $Z$ is superadditive. If $Cx_1 = z_1 \in Z(\beta_1)$ and $Cx_2 = z_2 \in Z(\beta_2)$, then $x_1 + x_2 \in \mathcal{X}(\beta_1 + \beta_2)$. Thus, there exists $z_3 \in Z(\beta_1 + \beta_2)$ such that $z_1 + z_2 \leqq z_3$. Also, there is no $z_4 \in Z(\beta_1 + \beta_2)$ such that $z_4 \leq z_1 + z_2$, because such a $z_4$ would not be a nondominated point for the MOIP with RHS $\beta_1 + \beta_2$.
\qed\endproof

Because the value function is set-valued, it is not immediate how to define an analog of single-objective superadditive duality. We present two variants, both of which coincide with the standard superadditive dual in the single-objective case.

\subsection{Set-valued Superadditive Dual} 
\label{sec:super_add_dual_set}

Recall from Section \ref{sec:background} that the superadditive IP dual contains the constraints $F(A_j) \geq c_j$. 
Our first formulation \eqref{eq:super_add_dual_set} generalizes this constraint to a set-valued counterpart. Consider the following problem.
\begin{equation}
\label{eq:super_add_dual_set}\tag{SDP}
 \begin{aligned}
 \min ~ & F(b)\\
 \text{s.t. } & \{C_j\} \preceqq F(A_j) & \text{for } j =1, \ldots, n,\\
 & 0 \in F(0), \\
 & F:\reals^m\rightarrow\mathcal{E} & \text{superadditive and nondecreasing}.
 \end{aligned}
\end{equation}

Let $\mathcal{F}$ be the set of functions feasible to \eqref{eq:super_add_dual_set}. Then, $F(b) \in \mathcal{E}$ for all $F \in \mathcal{F}$. We interpret the objective of \eqref{eq:super_add_dual_set} as finding the elements of the collection $\{ F(b) \ \vert\ F \in \mathcal{F} \} \subseteq \mathcal{E}$ that are nondominated from below with respect to the set-ordering ``$\preceqq$''. This objective is well-defined because ``$\preceqq$'' defines a partial order on $\mathcal{E}$. Proposition \ref{prop:SDP1_weak} establishes weak duality for \eqref{eq:super_add_dual_set}.

\begin{proposition}[Weak Duality for \eqref{eq:super_add_dual_set}]
\label{prop:SDP1_weak}
If $x \in \mathcal{X}(b)$ and $F$ is feasible to \eqref{eq:super_add_dual_set}, then $\{Cx\} \preceqq F(b)$.
\end{proposition}

\proof 

Let $x \in \mathcal{X}(b)$. As $F$ is nondecreasing and $Ax \leqq b$, we have $F(Ax) \preceqq F(b)$. Now, consider the set 
\begin{equation*} 
\begin{aligned}
S = \bigg\{ \sum\limits_{j = 1}^n z_j x_j \Big| z_j \in F(A_j) \bigg\}.
\end{aligned} 
\end{equation*}
We claim that $S \preceqq F(b)$. Because $x_j$ is a nonnegative integer, we can express $z_j x_j$ as $\sum\limits_{i = 1}^{x_j} z_j$. (If $x_j = 0$, the sum is empty and equal to $0 \in \mathbb{R}^k$.) Then, 
\begin{equation*} \begin{aligned}
S = \bigg\{ \sum\limits_{j = 1}^n \sum\limits_{i = 1}^{x_j} z_j \Big| z_j \in F(A_j) \bigg\}
\subseteq \sum\limits_{j=1}^n \sum\limits_{i = 1}^{x_j} F(A_j) = \sum\limits_{j=1}^n x_j F(A_j).
\end{aligned} 
\end{equation*}

Note that $Ax = \sum\limits_{j=1}^n A_j x_j$. Then, the superadditivity of $F$ implies that 
$\sum\limits_{j=1}^n x_j F(A_j) \preceqq F(Ax)$.
By Lemma \ref{lem:order_subset}, $S \preceqq F(Ax)$. By transitivity of the order, we have $S \preceqq F(b)$. It follows by Lemma \ref{lem:order_max} that $\Max(S) \preceqq F(b)$.

Next, we will show that $\{Cx\} \preceqq \Max(S)$. Because $F(A_j) \suceqq \{C_j\}$ for all $j = 1, \ldots, n$, there is a $z_j \in F(A_j)$ such that $C_j \leqq z_j$. 
Then, $w = \sum\limits_{j=1}^n z_j x_j \in S$ is such that $Cx \leqq w$. 
Therefore, there exists $\tilde{w} \in \Max(S)$ such that $\tilde{w} \geqq w \geqq Cx$. 

Further, suppose that $z_j \in F(A_j)$ are such that 
$\sum_{j = 1}^{n} z_j x_j \leq \sum_{j = 1}^{n}C_jx_j$.
We have already shown that there exists $\tilde{w} \in S$ such that $Cx \leqq \tilde{w}$. Therefore, 
\begin{equation*} 
\begin{aligned}
\sum_{j = 1}^{n} z_j x_j \leq \sum_{j = 1}^{n}C_j x_j = Cx \leqq \tilde{w},
\end{aligned} 
\end{equation*}
so that $\sum_{j = 1}^{n} z_j x_j \not \in \Max(S)$.
Thus, $\{Cx\} \preceqq \Max(S) \preceqq F(b)$, which implies $\{Cx\} \preceqq F(b)$.
\qed \endproof

We say that the dual problem \eqref{eq:super_add_dual_set} is  strong at an efficient solution $x^*$ of \eqref{eq:MOIP} if there is a function $F$ feasible to \eqref{eq:super_add_dual_set} such that $Cx^* \in F(b)$.
Theorems \ref{thm:SDP1_strong_supported} and \ref{thm:SDP1_single_strong_func} guarantee that \eqref{eq:super_add_dual_set} is strong at the supported primal solutions.

\begin{theorem}[Strong Duality for \eqref{eq:super_add_dual_set}]
\label{thm:SDP1_strong_supported}
If $x^*$ is a supported efficient solution for \eqref{eq:MOIP}, then there is a function $F$ feasible to \eqref{eq:super_add_dual_set} such that $Cx^* \in F(b)$. 
\end{theorem}

\proof 
Let $x^*$ be a supported efficient solution to \eqref{eq:MOIP} and let $\mu \in \reals^k_{>}$ be a supporting vector such that $x^*$ is optimal to the scalarized problem 
\begin{equation*}
\begin{aligned}
\max & \enskip \mu^\top Cx \quad
\text{s.t. } \enskip Ax \leqq b,\ x \in \Z_{\geqq}^n.
\end{aligned}
\end{equation*}
Let $z^* = Cx^*$. We need to show that there is a function $F$ feasible for \eqref{eq:super_add_dual_set} such that $z^* \in F(b)$.
By strong superadditive duality for single-objective IPs (Theorem \ref{thm:single_obj_strong_super}), there exists a nondecreasing superadditive function $f_\mu: \mathbb{R}^m \rightarrow \reals$ such that $f_\mu(0) = 0$, $f_\mu(A_j) \geq \mu^\top C_j$ for all $j$, and 
$f_\mu(b) = \mu^\top Cx^* = \mu^\top z^*.$
We use $f_\mu$ to define a function $F_\mu: \reals^m \rightarrow \mathcal{E}$ as follows:
\begin{equation*} \begin{aligned}
F_\mu(v) =
\{w\in \mathbb{R}^k \ \vert \ \mu^\top w = f_\mu(v)\}.
\end{aligned} \end{equation*}

We first show that $F_\mu$ is well-defined. Let $\hat{\mu} = \frac{1}{\|\mu\|^2}\mu$. Then, $F_\mu(v) \neq \emptyset$ for all $v$ because $(f_\mu(v) \hat{\mu}) \in F_\mu(v)$. Also, $0 \in F_\mu(0)$ because $f_\mu(0) = 0$. For arbitrary $v$, let $w_1,w_2 \in F_\mu(v)$ with $w_1 \not = w_2$. Then, $w_1 \not \leq w_2$ and $w_1 \not \geq w_2$ because $\mu^\top (w_1 - w_2) = 0$. Thus, $F_\mu(v) \in \mathcal{E}$.

Next, we show that $F_\mu$ is feasible to \eqref{eq:super_add_dual_set}.
Let $v_1, v_2 \in \reals^m$, $v_1 \leqq v_2$. Then, $f_\mu(v_1) \leqq f_\mu(v_2)$ as $f_\mu$ is nondecreasing. 
We need to show that $F_\mu(v_1) \preceqq F_\mu(v_2)$. 
Let $w_1 \in F_\mu(v_1)$.
Define $w_2 = w_1 + (f_\mu(v_2) - f_\mu(v_1))\hat{\mu}$ so that 
\[
\mu^\top w_2 = \mu^\top\left[w_1 {+} (f_\mu(v_2) - f_\mu(v_1))\hat{\mu}\right]= \mu^\top w_1 +(f_\mu(v_2) - \mu^\top w_1)\frac{\mu^\top \mu}{\|\mu\|^2} = f_\mu(v_2).
\]

\noindent Thus, $w_2 \in F_\mu(v_2)$. Because $\hat{\mu} \in \mathbb{R}^k_{>}$ and $(f_\mu(v_2) - f_\mu(v_1)) \geqq 0$, $w_1 \leqq w_2$.

Also, for any $w_3 \in F_\mu(v_2)$,
$\mu^\top w_1 =f_\mu(v_1) \leqq f_\mu(v_2) = \mu^\top w_3$, which implies that $w_1 \not \geq w_3$. Thus, $F_\mu(v_1) \preceqq F_\mu(v_2)$ and $F_\mu$ is nondecreasing.

Second, for any $v_1, v_2 \in \reals^m$, $f_\mu(v_1) + f_\mu(v_2) \leqq f_\mu(v_1 + v_2)$ because $f$ is superadditive. If $w_1 \in F_\mu(v_1), w_2 \in F_\mu(v_2)$, and $w_3 \in F_\mu(v_1 + v_2)$, then $\mu^\top w_3 \geqq \mu^\top (w_1 + w_2)$ so that $w_3 \not \leq w_1 + w_2$. Moreover, by choosing 
$
w_3 = (w_1 + w_2) + \left[f_\mu(v_1 + v_2) - f_\mu(v_1) - f_\mu(v_2)\right]\hat{\mu},
$
we have $w_3 \in F_\mu(v_1 + v_2)$ with $w_3 \geqq w_1 + w_2$. Thus, $F_\mu(v_1) + F_\mu(v_2) \preceqq F_\mu(v_1 + v_2)$, 
and $F_\mu$ is superadditive.

Finally, note that $f_\mu(A_j) \geqq \mu^\top C_j$. If $w \in F_\mu(A_j)$, then $\mu^\top w \geqq \mu^\top C_j$ so that $w \not \leq C_j$. Moreover, for
$
w = C_j + (f_\mu(A_j) - \mu^\top C_j)\hat{\mu},
$
we have $w \in F_\mu(A_j)$ and $w \geqq C_j$. Thus, $\{C_j\} \preceqq F_\mu(A_j)$.

Thus, $F_\mu$ is feasible to \eqref{eq:super_add_dual_set}. Because $\mu^\top z^* = f_\mu(b)$, $z^* \in F_\mu(b)$. 
\qed \endproof

Theorem \ref{thm:SDP1_strong_supported} does not address the behavior of \eqref{eq:super_add_dual_set} if the primal is infeasible. As in single-objective IPs \cite{wolsey2014integer}, the dual is unbounded in that case. The result is presented later in Corollary \ref{cor:infeasible_unbounded_set}, as the proof uses properties derived in Section \ref{sec:super_add_dual_others}. 

Theorem \ref{thm:SDP1_strong_supported} implies that for each primal supported solution $x$, there is a dual solution $G$ feasible to \eqref{eq:super_add_dual_set} such that $G(b)$ contains $Cx$. However, the theorem makes no statement about the number of feasible functions needed to obtain all such points. That is, Theorem \ref{thm:SDP1_strong_supported} alone does not guarantee the existence of a single function $F^*$ feasible to \eqref{eq:super_add_dual_set} such that $Cx^* \in F^*(b)$ for all supported efficient solutions $x^*$ to MOIP. Theorem \ref{thm:SDP1_single_strong_func} shows that such a function always exists. To prove this, we first establish that finitely many scalarizations suffice to recover all supported efficient solutions of \eqref{eq:MOIP}.

\begin{lemma}\label{lem:Finite_scalarizations}
There exist scalarizing vectors $\mu_1,\mu_2, \ldots, \mu_\ell \in \reals^k_{>}$ such that every supported efficient solution $x^*$ of \eqref{eq:MOIP} is an optimal solution of the scalarized problem 
\begin{equation}
\label{eq:scal_MOIP_y_i}\tag{MOIP$_{\mu_i}$}
\begin{aligned}
\max & \enskip \mu_i^\top Cx\quad
 \text{s.t. } \enskip Ax \leqq b,\  x \in \Z_+^n,
\end{aligned}
\end{equation}
for some $i = 1, \ldots, \ell$.
\end{lemma}

\proof
By Proposition \ref{prop:scal_convex}, the supported efficient solutions of \eqref{eq:MOIP} are precisely the integral efficient solutions to \eqref{eq:CH}. Because \eqref{eq:CH} is an MOLP, there are finitely many scalarizing vectors $\mu_1,\mu_2, \ldots, \mu_\ell \in \reals_{>}^k$ such that each efficient solution of \eqref{eq:CH} is an optimal solution to the LP
\begin{equation}
\label{eq:scal_MOLP_y_i}\tag{CH$_{y_i}$}
\begin{aligned}
\max & \enskip \mu_i^\top Cx\quad 
\text{s.t. } \enskip  x \in \mathrm{conv}(\mathcal{X}), 
\end{aligned}
\end{equation}
\noindent for some $i =1, \ldots, \ell$ \cite[Corollary 4.3.3]{luc_2016_Book}. So, if $x$ is a supported efficient solution of \eqref{eq:MOIP}, then $x$ is an integral optimal solution to \eqref{eq:scal_MOLP_y_i} for some $i =1, \ldots, \ell$ and therefore an optimal solution to \eqref{eq:scal_MOIP_y_i}. 
\qed \endproof

\begin{theorem}\label{thm:SDP1_single_strong_func}
There exists a function $F^*$ feasible to \eqref{eq:super_add_dual_set} such that if $x^*$ is a supported efficient solution of \eqref{eq:MOIP}, then $Cx^* \in F^*(b)$.
\end{theorem}

\proof 
Lemma \ref{lem:Finite_scalarizations} implies that there are finitely many vectors $\mu_1,\mu_2,\ldots, \mu_\ell \in \reals^k_{>}$ such that every supported efficient solution of \eqref{eq:MOIP} is an optimal solution to the scalarized IP \eqref{eq:scal_MOIP_y_i} for some $i = 1, \ldots, \ell$.

For every $i=1, \ldots \ell$, define $\hat{\mu_i} = \frac{\mu_i}{\|\mu_i\|^2}$. Then, there is a superadditive nondecreasing function $f_{\mu_i}:\reals^m \rightarrow \reals$ such that $f_{\mu_i}(A_j) \geq \mu_i^\top C_j$ for all $j = 1, \ldots, n$, and $\mu_i^\top Cx^* = f_{\mu_i}(b)$ for all $x^*$ optimal to \eqref{eq:scal_MOIP_y_i}. For $i = 1, \ldots, \ell$, consider the functions $F_{\mu_i}: \reals^m \rightarrow \mathcal{E}$ given by 
\begin{equation*} \begin{aligned}
F_{\mu_i}(v) =
\{w\in \mathbb{R}^k \ \vert \ \mu_i^\top w = f_{\mu_i}(v)\}
\end{aligned} \end{equation*}

\noindent 
As in the proof of Theorem \ref{thm:SDP1_strong_supported}, 
each $F_{\mu_i}(v)$ is feasible to \eqref{eq:super_add_dual_set}. Define the set $\mathcal{S}(v) = \cup_{i = 1}^{\ell}F_{\mu_i} (v)$ and consider the function 
\begin{equation*} 
\begin{aligned}
F^*(v) = \Min\ (\mathcal{S}(v)) = \Min\Big(\bigcup_{i = 1}^{\ell}F_{\mu_i} (v) \Big).
\end{aligned} 
\end{equation*}
\noindent 
We claim that $F^*$ is feasible to \eqref{eq:super_add_dual_set}.
For each $v$, the set $\mathcal{S}(v)$ is a finite union of hyperplanes in $\reals^k$. Therefore, either $F^*(v) = - M_{\infty}$ or $\mathcal{S}(v)$ has points that are nondominated from below. In either case, $F^*(v) \in \mathcal{E}$. 

Next, we show that $0 \in F^*(0)$. 
For each $i = 1,\ldots, \ell$, we have $0 \in F_{\mu_i}(0)$ and there is no $z \in F_{\mu_i}(0)$ with $z \leq 0$ because $F_{\mu_i}(0) \in \mathcal{E}$. In particular, there are no elements of $\mathcal{S}(0)$ which are dominated by $0$. Thus, $0$ is a nondominated point of $\mathcal{S}(0)$, that is, $0 \in F^*(0)$.

Because $F_{\mu_i}(A_j) \suceqq C_j$ for all $i =1, \ldots \ell$, there is no $z \in \mathcal{S}(A_j)$ such that $z \leq C_j$. In particular, because $F^*(A_j) \subseteq \mathcal{S}(A_j)$, there is no $z \in F^*(A_j)$ such that $z \leq C_j$. 

Next, we show that for every $j$, there is a $z \in F^*(A_j)$ such that $z \geqq C_j$.
Consider $\mathcal{Z}_j = \{z \in \mathcal{S}(A_j) \ \vert \ C_j \leqq z\}$. Note that $\mathcal{Z}_j \cap F_{\mu_i}(A_j) \neq \emptyset$ for all $i,j$ because each $F_{\mu_i}$ is feasible to \eqref{eq:super_add_dual_set}.
Let $z \in \Min(\mathcal{Z}_j)$, and we claim that $z \in F^*(A_j)$. Suppose if possible that there exists $w \in \mathcal{S}(A_j)$ such that $w \leq z$. 
Then there exists $i \in \{1, \ldots,\ell\}$ such that $w \in F_{\mu_i}(A_j)$. 
Then, $\tilde{z} = z+ (f_{\mu_i}(A_j) - \mu_i^\top z)\hat{\mu}_i$ satisfies $\tilde{z} \in F_{\mu_i}(A_j)$ and $z \leqq \tilde{z}$ because $z \in \Min (\mathcal{Z}_j)$ and $\tilde{z}$ is comparable with $z$. Then, the following three statements hold: $\mu_i^\top z \leqq \mu_i^\top \tilde{z}$, $\mu_i^\top w = \mu_i^\top \tilde{z}$, and $\mu_i^\top w < \mu_i^\top z$, which is a contradiction. Thus, there is no such $w$, and we have $z \in F^*(A_j)$. So, $\{C_j\} \preceqq F^*(A_j)$. 

To see that $F^*(v)$ is nondecreasing, let $v \leqq \tilde{v}$. Let $w \in F^*(v)$ and $\tilde{w} \in F^*(\tilde{v})$, and suppose without loss of generality that $\tilde{w} \in F_{\mu_1}(\tilde{v})$. Because $F_{\mu_1}$ is nondecreasing, $z \not \geq \tilde{w}$ for all $z \in F_{\mu_1}(v)$. If $w \in F_{\mu_i}(v)$ for $i \not = 1$, then 
\begin{equation*} \begin{aligned}
z = w + \left(f_{\mu_1}(v) - \mu_1^\top w\right)\hat{\mu_1}
\end{aligned} \end{equation*}
\noindent 
is an element of $F_{\mu_1}(v)$. If $w \in F^*(v)$, then $w \leqq z$ because $\hat{\mu_1}\in \reals^k_{>}$ and therefore $w$ and $z$ are comparable. It follows that $\tilde{w} \not \leq w$ because $w \leqq z$ and $\tilde{w} \not \leq z$ as $z \in F_{\mu_1}(v)$.

We now show that if $w \in F^*(v)$, then there is a $\tilde{w} \in F^*(\tilde{v})$ such that $\tilde{w} \geqq w$. For each $i = 1,\ldots,\ell$, define
\begin{equation*} \begin{aligned}
\tilde{w}_i = w + \left(f_{\mu_i}(\tilde{v}) -\mu_i^\top w\right)\hat{\mu_i},
\end{aligned} \end{equation*}
\noindent 
and note that $\tilde{w}_i \in F_{\mu_i}(\tilde{v})$. Moreover, because $F$ is nondecreasing and because $w \in F^*(v)$, we have $w \leqq \tilde{w}_i$ for all $i$. Let $\tilde{w} \in \Min(\{ \tilde{w}_i \ \vert\ i=1, \ldots,n\})$. Suppose there is some $i \in \{1,\ldots,\ell\}$ and $z \in F_{\mu_i}(\tilde{v})$ such that $z \leq \tilde{w}$. Then, if $\tilde{z} = \tilde{w} + \left(f_{\mu_i}(\tilde{v}) - \mu_i^\top \tilde{w}\right)\hat{\mu}_i$, it would follow that $\mu_i^\top z < \mu_i^\top \tilde{w}$, $\mu_i^\top z = \mu_i^\top \tilde{z}$, and $\mu_i \top \tilde{w} \leqq \mu_i^\top \tilde{z}$, which is a contradiction. So, no such $z$ exists and we have $\tilde{w} \in F^*(\tilde{v})$. Thus, $F^*(v) \preceqq F^*(\tilde{v})$. 

Next, we prove that $F^*$ is superadditive. Let $w_1 \in F^*(v_1)$ and $w_2 \in F^*(v_2)$, and $\tilde{w} \in F^*(v_1 + v_2)$. Suppose without loss of generality that $\tilde{w} \in F_{\mu_1}(v_1 + v_2)$. Then, for every $z_1 \in F_{\mu_1}(v_1)$ and $z_2 \in F_{\mu_1}(v_2)$, $\tilde{w} \not \leq z_1 + z_2$ by the superadditivity of $F_{\mu_1}$. In particular, setting
\begin{equation*} 
\begin{aligned}
z_1 = w_1 + \left(f_{\mu_1}(v_1) - \mu_{1}^\top w_1\right)\hat{\mu_1} 
\quad \text{ and } \quad 
z_2 = w_2 + \left(f_{\mu_1}(v_2) - \mu_1^\top w_2\right)\hat{\mu_1},
\end{aligned} 
\end{equation*}
we have $z_1 \in F_{\mu_1}(v_1)$, $z_2 \in F_{\mu_1}(v_2)$ and $w_1 + w_2 \leqq z_1 + z_2$. This in turn implies that $\tilde{w} \not \leq w_1 + w_2$ because $\tilde{w} \not \leq z_1 + z_2$.

Let $w_1 \in F^*(v_1), w_2 \in F^*(v_2)$, and let 
\begin{equation*} 
\begin{aligned}
\tilde{w} \in \Min\ \{(w_1 + w_2) + \left(f_{\mu_i}(v_1 + v_2) - \mu_i^\top(w_1 + w_2)\right)\hat{\mu_i} \ \vert \ 1 \leq i \leq \ell\}.
\end{aligned} 
\end{equation*}
Then, for each $i$, there is a $z_2 \in F_{\mu_i}(v_1 + v_2)$ such that $\tilde{w} \leqq z_2$. Therefore, there is no $z\in F_{\mu_i}(v_1 + v_2)$ with $z \leq \tilde{w}$ because if there were, then $\mu_i^\top z < \mu_i^\top \tilde{w}$, $\mu_i^\top \tilde{w} \leqq \mu_i^\top z_2$, and $\mu_i^\top z = \mu_i^\top z_2$, which cannot happen. So, $\tilde{w} \in F^*(v_1 + v_2)$. Therefore, $F^*(v_1) + F^*(v_2) \preceqq F^*(v_1 + v_2)$. 

So, $F^*$ is feasible to \eqref{eq:super_add_dual_set}. For each supported primal efficient solution $x^*$, there exists a $1 \leq i \leq \ell$ such that $Cx^* \in F_{\mu_i}(b)$. Because $\{Cx^*\} \preceqq F^*(b)$ by Proposition \ref{prop:SDP1_weak}, it follows that $Cx^* \in F^*(b)$.
\qed \endproof

Thus, the superadditive dual \eqref{eq:super_add_dual_set} is a strong dual to \eqref{eq:MOIP} at supported primal efficient solutions. However, the set-valued functions it employs are difficult to characterize, which limits the immediate algorithmic utility of this dual in providing a bound set for search-based solution methods. Several researchers have developed methods to approximate the single-objective IP value function, and future research in this area could extend those methods to approximate the set-valued MOIP value function. 
In this paper, we consider a restricted dual that only includes vector-valued functions.


\subsection{Vector-Valued Superadditive Dual}
\label{sec:super_add_dual_others}

We formulate another dual problem to \eqref{eq:MOIP} by restricting the feasible region in \eqref{eq:super_add_dual_set} to include only vector-valued functions. In other words, we consider only those functions $F$ feasible to \eqref{eq:super_add_dual_set} for which $F(v)$ is a singleton set in $\reals^k$ for all $v \in \reals^m$. In this case, we denote $F$ as a $k$-tuple with components $f_i$, $i = 1, \ldots, k$, where $f_i : \reals^m \rightarrow \reals$, and have the following dual formulation

\begin{equation}
\label{eq:super_add_dual_vec}\tag{VSDP}
 \begin{aligned}
 \min \enskip & F(b) = (f_1(b), \ldots, f_k(b))^\top \\
 \text{s.t.} \enskip & f_i(A_j) \geqq c_{ij} & \quad \text{for all } i,j,\\
 & f_i(0) = 0 & \quad \text{for all } i,\\
 & f_i:\mathbb{R}^m\rightarrow \mathbb{R} &\quad \text{nondecreasing and superadditive for all $i$.}
 \end{aligned}
\end{equation}

\begin{proposition}[Weak Duality for \eqref{eq:super_add_dual_vec}]\label{prop_super_weak_dual}
If $F$ is feasible to \eqref{eq:super_add_dual_vec}, then $Cx \leqq F(b)$ for all $x \in \mathcal{X}(b)$. 
\end{proposition}

\proof 
For any $x \in \mathcal{X}(b)$,
\begin{equation*} \begin{aligned}
\begin{aligned}
Cx = \begin{pmatrix}
\sum_{j = 1}^{n}c_{1,j}x_j\\
\vdots\\
\sum_{j = 1}^{n}c_{k,j}x_j\\
\end{pmatrix}
\leqq \begin{pmatrix}
\sum_{j = 1}^{n}f_1(A_j)x_j\\
\vdots\\
\sum_{j = 1}^{n}f_k(A_j)x_j\\
\end{pmatrix}
\leqq \begin{pmatrix}
f_1(Ax)\\
\vdots\\
f_k(Ax)
\end{pmatrix}
\leqq F(b),
\end{aligned}
\end{aligned} \end{equation*}
where the first inequality follows from the $f_i(A_j) \geqq c_{i,j}$ constraint, the second inequality is due to $f_i$ being superadditive, and the third inequality follows from $f_i$ being nondecreasing.
\qed \endproof

As in \eqref{eq:super_add_dual_set}, the minimization in \eqref{eq:super_add_dual_vec} is with respect to the partial order $\preceqq$. However, if $\{F_\alpha(b) \ \vert\ \alpha \in \mathcal{I} \}$ is the set of nondominated points of \eqref{eq:super_add_dual_vec} for some indexing set $\mathcal{I}$, then $F_{\alpha_1}(b)$ and $F_{\alpha_2}(b)$ are incomparable under $\preceqq$ for all $\alpha_1,\alpha_2 \in \mathcal{I}$, $\alpha_1 \neq \alpha_2$. Because $F_{\alpha}(b)$ is a singleton set for each $\alpha \in \mathcal{I}$, we can identify $F_\alpha(b)$ with its sole element $z_\alpha$. Under this identification, the nondominated set $\{F_\alpha(b) \ \vert \ \alpha \in \mathcal{I} \}$ is equivalent to $\{z_\alpha \ \vert\ \alpha \in \mathcal{I}\} \in \mathcal{E}$. Moreover,
\[
\{z_\alpha \ \vert \ \alpha \in \mathcal{I}\} = \Min\Big( \underset{F \text{ feasible to \eqref{eq:super_add_dual_vec}} }{\bigcup} F(b)\Big). 
\]
In this way, \eqref{eq:super_add_dual_vec} is equivalent to solving a multiobjective problem whose objective values are elements of $\reals^k$. For the remainder of this section, we view \eqref{eq:super_add_dual_vec} as a problem in $\reals^k$ and make the identification of $F(b)$ with its sole element for a singleton set $F(b)$. 

If an objective of \eqref{eq:MOIP} is unbounded, then the corresponding objective of \eqref{eq:super_add_dual_vec} is infeasible by Proposition \ref{prop:IP_weak_SA}, so that \eqref{eq:super_add_dual_vec} is infeasible. 
Proposition \ref{prop:vec_dual_ideal} shows that the upper bound due to \eqref{eq:super_add_dual_vec} is tighter than any other singleton upper bound set for the MOIP. 

\begin{proposition}\label{prop:vec_dual_ideal}
Let $y \in \mathbb{R}^k$ such that $Cx \leqq y$ for all $x \in \mathcal{X}(b)$. Then, there exists a feasible solution $F$ of \eqref{eq:super_add_dual_vec} such that $F(b) \leqq y$. Moreover, if \eqref{eq:super_add_dual_vec} has efficient solutions, then there is an efficient solution $F^*$ of \eqref{eq:super_add_dual_vec} such that $F^*(b) \leqq y$. 
\end{proposition}

\proof 
If $Cx \leqq y$ for all $x \in \mathcal{X}(b)$, then for each $1\leq i\leq k$, $c_ix \leq y_i$ is a valid inequality. This implies that there is a superadditive nondecreasing function $f_i:\mathbb{R}^m \rightarrow \reals$ such that for each $1\leq j\leq n$, $f_i(A_j) \geq c_{i,j}$, $f_i(0) = 0$, and $f_i(b) \leq y_i$. Then, $F = (f_1, \ldots, f_k)^\top$ is feasible to \eqref{eq:super_add_dual_vec}. So, if \eqref{eq:super_add_dual_vec} has efficient solutions, then there exists $F^*$ efficient to \eqref{eq:super_add_dual_vec} with
$F^*(b) \leqq F(b) \leqq y$.
\qed \endproof

Proposition \ref{prop:vec_dual_ideal} has several corollaries that describe the nondominated points of \eqref{eq:super_add_dual_vec}. 

\begin{corollary}\label{cor:infeasible_unbounded_vec}
If \eqref{eq:MOIP} is infeasible, then \eqref{eq:super_add_dual_vec} is unbounded.
\end{corollary}

\proof
If \eqref{eq:MOIP} is infeasible, then for each $y \in \reals^k$, $Cx \leqq y$ for all $x \in \mathcal{X}(b)$. So, by Proposition \ref{prop:vec_dual_ideal}, for each $y \in \reals^k$ there is a feasible objective value $F(b)$ to \eqref{eq:super_add_dual_set} with $F(b) \leqq y$.
\qed\endproof

\begin{corollary}\label{cor:infeasible_unbounded_set}
If \eqref{eq:MOIP} is infeasible, then \eqref{eq:super_add_dual_set} is unbounded.
\end{corollary}

\proof
\eqref{eq:super_add_dual_set} is a relaxation of \eqref{eq:super_add_dual_vec}. Corollary \ref{cor:infeasible_unbounded_vec} therefore implies that if \eqref{eq:MOIP} is infeasible, then for each $y \in \reals^k$ there is a feasible objective value $F(b)$ to \eqref{eq:super_add_dual_set} with $F(b) \preceqq \{y\}$.
\qed \endproof

\begin{corollary}\label{cor:unique_nondom}
If $G^*(b)$ is a nondominated point of \eqref{eq:super_add_dual_vec}, then $G^*(b)$ is the unique nondominated point of \eqref{eq:super_add_dual_vec}.
\end{corollary}

\proof 
Let $G^*(b)$ be a nondominated point of \eqref{eq:super_add_dual_vec}. This implies that \eqref{eq:super_add_dual_vec} is feasible and each single-objective IP $\max~\{c_ix \ \vert \ x \in \mathcal{X}(b)\}$ has a finite optimal value. Let $y \in \reals^k$ be the vector with components $y_i = \max~\{c_ix \ \vert \ x \in \mathcal{X}(b)\}$. Proposition \ref{prop_super_weak_dual} then implies that for each $i \in \{1,2,\ldots, k\}$, the $i^{th}$ objective of $G^*(b)$ is bounded below by $y_i$. That is,  $y \leqq G^*(b)$.

On the other hand, $Cx \leqq y$ holds for all $x \in \mathcal{X}(b)$ so that Proposition \ref{prop:vec_dual_ideal} implies the existence of a feasible function $F^*(b)$ to \eqref{eq:super_add_dual_vec} such that $F^*(b) \leqq y$. Then, $F^*(b) \leqq G^*(b)$ by the transitivity of $\leqq$. Because $F^*(b)$ and $G^*(b)$ are both nondominated points of \eqref{eq:super_add_dual_vec}, this inequality must hold with equality. Because $G^*(b)$ was an arbitrary nondominated point, this implies that \eqref{eq:super_add_dual_vec} has a unique nondominated point. 
\qed \endproof

\begin{remark}\label{rem:ideal}
    Proposition \ref{prop:vec_dual_ideal} and Corollary \ref{cor:unique_nondom} imply that \eqref{eq:super_add_dual_vec} computes the ideal point $y^I$ of \eqref{eq:MOIP}.
\end{remark}

One approach to proving strong duality of \eqref{eq:super_add_dual_vec} is to show that Proposition \ref{prop:vec_dual_ideal} holds for every $z \in Z(b)$. This, however, may not be true as an element of $Z(b)$ may be incomparable with some $Cx$ in the value-set of the MOIP. This is illustrated in Example \ref{ex_MOIP_val_func}. 

\begin{example}\label{ex_MOIP_val_func}
Consider again the MOIP from Example \ref{ex:P_not_open}.
\begin{equation*} \begin{aligned}
\max & \enskip \begin{bmatrix} 1 & -\frac{1}{2}\\ -\frac{1}{2} & 1\end{bmatrix} \begin{bmatrix} x_1\\x_2 \end{bmatrix}\quad 
\text{s.t.} \enskip x_1 + x_2 \leqq 1,\ x_1,x_2 \in \mathbb{Z}_{+},
\end{aligned} \end{equation*}
\noindent 
This problem has nondominated points $\left\{(1,-0.5)^\top, (0,0)^\top, (-0.5,1)^\top\right\} = Z(1)$. But for any $z \in Z(1)$, there is a feasible $x$ such that $Cx \not \leqq z$, where $C$ is the objective matrix. For example, if $z = (1, -0.5)^\top$, then $x = (0,0)^\top$ is feasible but $Cx \not \leqq z$. So, there is no $z\in Z(1)$ such that $Cx \leqq z$ holds for all $x \in \mathcal{X}$.
\qed\end{example}

In the special case when \eqref{eq:MOIP} has a unique nondominated point, all elements of the set of feasible objective values must be comparable with the sole element in the singleton set $Z(b)$. This is established in Lemma \ref{lemma_ideal_valid}. We slightly abuse the notation $Z(b)$ to also denote the unique nondominated point contained in the set $Z(b)$.

\begin{lemma}\label{lemma_ideal_valid}
If \eqref{eq:MOIP} has a unique nondominated point, then $Cx \leqq Z(b)$ for all $x \in \mathcal{X}(b)$. 
\end{lemma}

\proof 
Let $x^*$ be an efficient solution of \eqref{eq:MOIP}, so that $Cx^* = Z(b)$. Suppose for the sake of a contradiction that there was an $x \in \mathbb{Z}^n_{+}$ such that $Ax \leqq b$ but $Cx \not \leqq Z(b)$. Then, there would be an objective $i \in \{1, \ldots, k\}$ such that $c_ix > c_ix^*$, which contradicts the hypothesis that $Cx^*$ is the unique nondominated point.
\qed \endproof

The statement $Cx \leqq Z(b)$ may be ill-defined if \eqref{eq:MOIP} does not have a unique nondominated point. So, the direct converse of Lemma \ref{lemma_ideal_valid} is not well defined, but a modified converse holds.

\begin{lemma}\label{lem:ideal_if_valid}
If there exists $z \in Z(b)$ such that $Cx \leqq z$ for all $x \in \mathcal{X}(b)$, then \eqref{eq:MOIP} has a unique nondominated point.
\end{lemma}

\proof 
Let $z_1,z_2 \in Z(b)$ such that $Cx \leqq z_1$ is a valid inequality. Because $z_1,z_2$ are feasible objective values, there are $x_1,x_2 \in \mathcal{X}(b)$ such that $z_1 = Cx_1$ and $z_2 = Cx_2$. If $Cx \leqq z_1$ is a valid inequality, then $z_1 \geqq Cx_2 = z_2$. On the other hand, if $z_2 \in Z(b)$, then $z_1 \not \geq z_2$, which implies that $z_1 = z_2$. Because $z_2 \in Z(b)$ was arbitrary, this implies that $Z(b)$ has only one element. Therefore, \eqref{eq:MOIP} has a unique nondominated point.
\qed \endproof

Theorem \ref{thm_ideal_super_iff} completely characterizes the strength of \eqref{eq:super_add_dual_vec}.

\begin{theorem}\label{thm_ideal_super_iff}

The dual problem \eqref{eq:super_add_dual_vec} is strong if and only if the primal problem \eqref{eq:MOIP} has a unique nondominated point. 
\end{theorem}

\proof 

Suppose that \eqref{eq:super_add_dual_vec} is strong. Then, there exists an efficient solution $F^*$ to \eqref{eq:super_add_dual_vec} such that $F^*(b) \in Z(b)$. Then, Proposition \ref{prop_super_weak_dual} implies that $Cx \leqq F^*(b)$ is a valid inequality. Because $F^*(b) \in Z(b)$, Lemma \ref{lem:ideal_if_valid} implies that \eqref{eq:MOIP} has a unique nondominated point. 

Conversely suppose that \eqref{eq:MOIP} has a unique nondominated point $Z(b) = \left( z_1(b), \ldots, z_k(b) \right)^\top$. Then, $Cx \leqq Z(b)$ is a valid inequality by Lemma \ref{lemma_ideal_valid}. Proposition \ref{prop:vec_dual_ideal} implies that there is a feasible solution $F(b)$ to \eqref{eq:super_add_dual_vec} with $F(b) \leqq Z(b)$. Proposition \ref{prop_super_weak_dual} implies that this inequality must hold with equality. Moreover, $F(b)$ must be an efficient solution because otherwise there would be a feasible solution $G(b)$ with $G(b) \leq F(b) \leqq Z(b)$, a contradiction with Proposition \ref{prop_super_weak_dual}. Then, $F(b)$ is an efficient solution to \eqref{eq:super_add_dual_vec} with $F(b) = Z(b)$. 
\qed \endproof

In the single-objective case, every \eqref{eq:MOIP} has a unique nondominated point, so that Theorem \ref{thm_ideal_super_iff} coincides with Theorem \ref{thm:single_obj_strong_super}. 

\subsubsection{MOLP reformulation of \eqref{eq:super_add_dual_vec}}
The superadditive dual of a single-objective IP can be formulated as an (exponentially large) LP in the special case where all entries of $A$ and $b$ are nonnegative integers. The subsequent discussion shows that \eqref{eq:super_add_dual_vec} can be cast as an MOLP in a similar manner. 

If \eqref{eq:super_add_dual_vec} has a nondominated point, it must be unique. Then, there exists a dual efficient function $F^* = (f_1^*, \ldots, f_k^*)^\top$ such that each $f_i^*$ is an optimal solution to the problem 
\begin{equation*} 
\begin{aligned}
\min \{f_i(b)\ \vert \ f_i(0) = 0, f_i(A_j) \geq c_{i,j}, 
f_i \text{ superadditive and nondecreasing}\},
\end{aligned} 
\end{equation*}
which is the superadditive dual for the single-objective IP $\max\{c_ix \ \vert \ x \in \mathcal{X}(b)\}$. If $A$ and $b$ have all nonnegative integral entries, then $f_i^*(b)$ is the optimal value of the following LP (see \cite{wolsey1981integer} for details).
\begin{equation*} 
\begin{aligned}
 \min \enskip &f_i(b) \\ 
 \text{s.t. } \enskip & f_{i}(A_j) \geqq c_{i,j} & 1 \leq j \leq n,\\
 & f_{i}(d_1) + f_i(d_2) - f_i(d_1 + d_2) \leqq 0 & \enskip \text{for all } 0\leqq d_1, d_2, (d_1 + d_2) \leqq b,\\
 & f_i(0) = 0, f_i(d) \geqq 0.
\end{aligned} 
\end{equation*}

This leads to the following MOLP reformulation of \eqref{eq:super_add_dual_vec} for instances in which all entries of $A$ and $b$ are nonnegative integers.
\begin{equation}
\label{eq:SDLP}\tag{SDMOLP}
\begin{aligned}
 \min\enskip & F(b) = \left(f_1(b), \ldots, f_k(b)\right)^\top\\
 \text{s.t. } \enskip & f_{i}(A_j) \geqq c_{i,j} & \text{for all } i, j,\\
 & f_{i}(d_1) + f_i(d_2) - f_i(d_1 + d_2) \leqq 0 & \enskip \text{for all } 0\leqq d_1, d_2, (d_1 + d_2) \leqq b,\\
 & f_i(0) = 0, f_i(d) \geqq 0.
\end{aligned}
\end{equation}

Because \eqref{eq:super_add_dual_vec} (and therefore \eqref{eq:SDLP}) has a unique nondominated point, a single scalarization suffices to recover an efficient solution to \eqref{eq:SDLP}. Thus, the vector-valued dual of this problem is an MOLP with a unique nondominated point. However, this dual is not strong unless the MOIP itself has a unique nondominated point. 

\begin{remark}
From Remark \ref{rem:ideal} and the above discussion, it follows that the ideal point $y^I$ for an MOIP with nonnegative constraint and objective coefficients can be obtained by solving a single (large) LP. 
\end{remark}

Example \ref{ex:knapsack} illustrates the MOLP reformulation on a bi-objective knapsack problem. Note that this MOLP has a unique nondominated point. Because the primal problem does not have a unique nondominated point, the vector-valued dual is not strong. Nonetheless, its (unique) nondominated point provides an upper bound on $\Max(\mathcal{Y}_{\mathrm{MOIP}})$. 

\begin{example}
\label{ex:knapsack}
Consider the MOIP
\begin{equation*} 
\begin{aligned}
\max \enskip & \begin{bmatrix} 2 & 1\\ 1 & 2\end{bmatrix}\begin{bmatrix} x_1\\ x_2 \end{bmatrix}\quad 
\text{s.t. } \enskip x_1 + x_2 \leqq 2,\  x_1,x_2 \in \mathbb{Z}_{+}.
\end{aligned}
\end{equation*}
\noindent
Following the steps described above, the MOLP formulation of the vector-valued superadditive dual for this problem is 
\begin{equation*} 
\begin{aligned}
\min \enskip & F(2) = (f_1(2), f_2(2))^\top \\
\text{s.t. } \enskip & \left(f_1(1), f_2(1)\right)^\top &\geqq (2,1)^\top,\\
& (f_1(1),f_2(1))^\top &\geqq (1,2)^\top,\\
& 2(f_1(1),f_2(1))^\top - (f_1(2),f_2(2))^\top &\leqq (0,0)^\top.
\end{aligned} 
\end{equation*}

\noindent This MOLP has a unique nondominated point $F^*(2) = (4, 4)^\top$. On the other hand, the original MOIP has nondominated points $\left\{(4,2)^\top, (3,3)^\top, (2,4)^\top\right\}$. 
\qed \end{example}

\section{Conclusion}\label{sec:conclusion}
In this paper, we analyzed relaxations and developed a duality framework for MOIPs by leveraging results from single-objective integer programming. In particular, we formulated the Lagrangian relaxation of an MOIP and compared it with the continuous and convex hull relaxations. The convex hull relaxation is tight at supported efficient solutions of the MOIP but not at unsupported solutions. We showed via an example that a Lagrangian relaxation can provide a tighter upper bound at unsupported nondominated points. 

We presented an MOIP Lagrangian dual that generalizes the single-objective counterpart, relying on the idea of finding the best upper bound over all Lagrangian relaxations. The behavior of this dual at supported solutions, including conditions for strong duality, mimics those derived in the single-objective case; the analysis is aided by scalarization techniques. The properties of the dual are harder to analyze at unsupported solutions. This is due in part to the complicated geometry of the dual feasible set $\mathcal{Y}_{\mathrm{LD}}$. Every point in the primal nondominated set has an upper bound in the dual feasible set $\mathcal{Y}_{\mathrm{LD}}$, but the non-convexity of $\mathcal{Y}_{\mathrm{LD}}$ implies that it also contains elements that are incomparable with the primal nondominated points. 
In particular, $\Min(\mathcal{Y}_{\mathrm{LD}})$ does not necessarily provide the ``best'' upper bound on the primal nondominated points and the Lagrangian relaxations themselves may be more informative in this respect. 
We presented computational evidence to illustrate that a na\"ive approximation to the Lagrangian dual bound set can provide a tighter upper bound than one obtained via convex hull relaxation. 

We also introduced two superadditive duals, namely, a set-valued formulation and a vector-valued variant. The set-valued problem considers set-valued functions that are non-decreasing and superadditive, inspired by the properties of the MOIP value function. This dual is strong at supported efficient solutions of the primal. 
The vector-valued dual is constructed by restricting the set-valued dual to functions from $\reals^m$ to $\reals^k$; it is strong if and only if the primal has a unique nondominated point. Given any upper bound $z$ on the set of feasible primal objective values, there exists a vector-valued dual feasible solution that provides a tighter upper bound. 
In the special case where the constraint parameters are nonnegative integers, the vector-valued dual can be formulated as an MOLP. Notably, the vector-valued superadditive dual provides an alternate method for computing the ideal point of an MOIP via a single (large) LP in case of nonnegative problem parameters. 


Our computational experiments have promising results, but our approach of enumerating a set of Lagrange multipliers over an equispaced grid does not scale well to larger problems. Future work in this area could focus on algorithmic aspects of the Lagrangian dual, especially with a view on techniques for selecting Lagrange multipliers. The IP value function is hard to compute even for the single-objective IP in the general case, but several researchers have developed algorithms for value function approximation for structured IPs. Extension of these methods to the multiobjective superadditive dual offers another promising avenue for future research.


%
%

\bibliographystyle{spmpsci} 
\bibliography{ref.bib} 

\begin{thebibliography}{10}
\providecommand{\url}[1]{{#1}}
\providecommand{\urlprefix}{URL }
\expandafter\ifx\csname urlstyle\endcsname\relax
  \providecommand{\doi}[1]{DOI~\discretionary{}{}{}#1}\else
  \providecommand{\doi}{DOI~\discretionary{}{}{}\begingroup
  \urlstyle{rm}\Url}\fi

\bibitem{AnejaNair79BicriteriaTransport}
Aneja, Y.P., Nair, K.P.K.: Bicriteria transportation problem.
\newblock Management Science \textbf{25}(1), 73--78 (1979)

\bibitem{Benson2009}
Benson, H.P.: Multi-objective optimization: {P}areto optimal solutions,
  properties.
\newblock In: C.A. Floudas, P.M. Pardalos (eds.) Encyclopedia of Optimization,
  pp. 2478--2481. Springer US, Boston, MA (2009)

\bibitem{BOLAND2017OptimizeoverEfficient}
Boland, N., Charkhgard, H., Savelsbergh, M.: A new method for optimizing a
  linear function over the efficient set of a multiobjective integer program.
\newblock European Journal of Operational Research \textbf{260}(3), 904--919
  (2017)

\bibitem{CERQUEUS2015SurrogateBoundSets}
Cerqueus, A., Przybylski, A., Gandibleux, X.: Surrogate upper bound sets for
  bi-objective bi-dimensional binary knapsack problems.
\newblock European Journal of Operational Research \textbf{244}(2), 417--433
  (2015)

\bibitem{CORLEY1984MATRIXDUAL}
Corley, H.: Duality theory for the matrix linear programming problem.
\newblock Journal of Mathematical Analysis and Applications \textbf{104}(1),
  47--52 (1984)

\bibitem{DACHERT2017EfficientComputationSearch}
Dächert, K., Klamroth, K., Lacour, R., Vanderpooten, D.: Efficient computation
  of the search region in multi-objective optimization.
\newblock European Journal of Operational Research \textbf{260}(3), 841--855
  (2017)

\bibitem{Ehrgott2005MultiCriteria}
Ehrgott, M.: Multicriteria Optimization.
\newblock Springer Berlin, Heidelberg (2005)

\bibitem{Ehrgott_Scalarization}
Ehrgott, M.: A discussion of scalarization techniques for multiple objective
  integer programming.
\newblock Annals of Operations Research \textbf{147}(1), 343--360 (2006)

\bibitem{ehrgott2001bounds}
Ehrgott, M., Gandibleux, X.: Bounds and bound sets for biobjective
  combinatorial optimization problems.
\newblock In: Multiple Criteria Decision Making in the New Millennium, pp.
  241--253. Springer (2001)

\bibitem{EHRGOTT2007BoundSets}
Ehrgott, M., Gandibleux, X.: Bound sets for biobjective combinatorial
  optimization problems.
\newblock Computers \& Operations Research \textbf{34}(9), 2674--2694 (2007)

\bibitem{Fisher1981Lagrange}
Fisher, M.L.: The {L}agrangian relaxation method for solving integer
  programming problems.
\newblock Management Science \textbf{27}(1), 1--18 (1981)

\bibitem{FORGET2022909}
Forget, N., Gadegaard, S.L., Nielsen, L.R.: Warm-starting lower bound set
  computations for branch-and-bound algorithms for multiobjective integer
  linear programs.
\newblock European Journal of Operational Research \textbf{302}(3), 909--924
  (2022)

\bibitem{Gale1951LinearProgrammingGames}
Gale, D., Kuhn, H.W., Tucker, A.W.: Linear programming and the theory of games.
\newblock Activity Analysis of Production and Allocation \textbf{13}, 317--335
  (1951)

\bibitem{vOptSolver2021}
Gandibleux, X., Soleihac, G., Przybylski, A.: v{O}pt{S}olver: an ecosystem for
  multi-objective linear optimization.
\newblock In: JuliaCon 2021 (2021)

\bibitem{vOptSolver2017MCDM}
Gandibleux, X., Soleilhac, G., Przybylski, A., Lucas, F., Ruzika, S.,
  Halffmann, P.: v{O}pt{S}olver, a ``get and run" solver of multiobjective
  linear optimization problems built on {J}ulia and {JuMP}.
\newblock In: MCDM2017: 24th International Conference on Multiple Criteria
  Decision Making, vol.~88 (2017)

\bibitem{vOptSolver2017IFORS}
Gandibleux, X., Soleilhac, G., Przybylski, A., Ruzika, S.: v{O}pt{S}olver: an
  open source software environment for multiobjective mathematical
  optimization.
\newblock In: IFORS2017: 21st Conference of the International Federation of
  Oprational Research Societies (2017)

\bibitem{GEOFFRION1968ProperEfficiency}
Geoffrion, A.M.: Proper efficiency and the theory of vector maximization.
\newblock Journal of Mathematical Analysis and Applications \textbf{22}(3),
  618--630 (1968)

\bibitem{Geoffrion1974Lagrange}
Geoffrion, A.M.: Lagrangean relaxation and its uses in integer programming.
\newblock Mathematical Programming \textbf{2}, 82--114 (1974)

\bibitem{GurionLuc2014MOLPLagrange}
Gourion, D., Luc, D.: Saddle points and scalarizing sets in multiple objective
  linear programming.
\newblock Mathematical Methods of Operations Research \textbf{80}(1), 1--27
  (2014)

\bibitem{Haimes1971epsilon_constraint}
Haimes, Y., Lasdon, L., Wismer, D.: On a bicriterion formulation of the
  problems of integrated system identification and system optimization.
\newblock IEEE Transactions on Systems, Man, and Cybernetics \textbf{SMC-1}(3),
  296--297 (1971)

\bibitem{Halffmann2022Review}
Halffmann, P., Schäfer, L.E., Dächert, K., Klamroth, K., Ruzika, S.: Exact
  algorithms for multiobjective linear optimization problems with integer
  variables: A state of the art survey.
\newblock Journal of Multi-Criteria Decision Analysis pp. 1--23 (2022)

\bibitem{Hamel2004MOLPLagrange}
Hamel, A.H., Heyde, F., L{\"o}hne, A., Tammer, C., Winkler, K.: Closing the
  duality gap in linear vector optimization.
\newblock Journal of Convex Analysis \textbf{11}(1), 163--178 (2004)

\bibitem{Heyde2008MOLPGeoDual}
Heyde, F., L\"ohne, A.: Geometric duality in multiple objective linear
  programming.
\newblock SIAM Journal on Optimization \textbf{19}, 836--845 (2008)

\bibitem{Heyde2009SetDual}
Heyde, F., Löhne, A., Tammer, C.: Set-valued duality theory for multiple
  objective linear programs and application to mathematical finance.
\newblock Mathematical Methods of Operations Research \textbf{69}(1), 159--179
  (2009)

\bibitem{Hooker2009}
Hooker, J.N.: Integer programming duality.
\newblock In: C.A. Floudas, P.M. Pardalos (eds.) Encyclopedia of Optimization,
  pp. 1657--1667. Springer US, Boston, MA (2009)

\bibitem{Iserman74ProperEfficiency}
Isermann, H.: Proper efficiency and the linear vector maximum problem.
\newblock Operations Research \textbf{22}(1), 189--191 (1974)

\bibitem{isermann1978some}
Isermann, H.: On some relations between a dual pair of multiple objective
  linear programs.
\newblock Zeitschrift f{\"u}r Operations Research \textbf{22}(1), 33--41 (1978)

\bibitem{JEROSLOW1978121}
Jeroslow, R.: Cutting-plane theory: Algebraic methods.
\newblock Discrete Mathematics \textbf{23}(2), 121--150 (1978)

\bibitem{jozefowiez2012generic}
Jozefowiez, N., Laporte, G., Semet, F.: A generic branch-and-cut algorithm for
  multiobjective optimization problems: Application to the multilabel traveling
  salesman problem.
\newblock INFORMS Journal on Computing \textbf{24}(4), 554--564 (2012)

\bibitem{Klamroth_dual}
Klamroth, K., Tind, J., Zust, S.: Integer programming duality in multiple
  objective programming.
\newblock Journal of Global Optimization \textbf{29}(1), 1--18 (2004)

\bibitem{kornbluth1974duality}
Kornbluth, J.: Duality, indifference and sensitivity analysis in multiple
  objective linear programming.
\newblock Journal of the Operational Research Society \textbf{25}(4), 599--614
  (1974)

\bibitem{Lohne2011}
L\"ohne, A.: Vector Optimization with Infimum and Supremum.
\newblock Springer-Verlag, Berlin (2011)

\bibitem{Luc2011MOLPDual}
Luc, D.T.: On duality in multiple objective linear programming.
\newblock European Journal of Operational Research \textbf{210}(2), 158--168
  (2011)

\bibitem{luc_2016_Book}
Luc, D.T.: Multiobjective Linear Programming: An Introduction.
\newblock Springer (2016)

\bibitem{lust2010twophasePareto}
Lust, T., Teghem, J.: Two-phase {P}areto local search for the biobjective
  traveling salesman problem.
\newblock Journal of Heuristics \textbf{16}(3), 475--510 (2010)

\bibitem{machuca2016lower}
Machuca, E., Mandow, L.: Lower bound sets for biobjective shortest path
  problems.
\newblock Journal of Global Optimization \textbf{64}(1), 63--77 (2016)

\bibitem{glpk}
Makhorin, A.: {GLPK} ({GNU} linear programming kit).
\newblock \url{https://www.gnu.org/software/glpk}  (2012)

\bibitem{mavrotas2005multi}
Mavrotas, G., Diakoulaki, D.: Multi-criteria branch and bound: A vector
  maximization algorithm for mixed 0-1 multiple objective linear programming.
\newblock Applied Mathematics and Computation \textbf{171}(1), 53--71 (2005)

\bibitem{ozpeynirci2010exact}
{\"O}zpeynirci, {\"O}., K{\"o}ksalan, M.: An exact algorithm for finding
  extreme supported nondominated points of multiobjective mixed integer
  programs.
\newblock Management Science \textbf{56}(12), 2302--2315 (2010)

\bibitem{PRZYBYLSKI2017MOBranchBound}
Przybylski, A., Gandibleux, X.: Multi-objective branch and bound.
\newblock European Journal of Operational Research \textbf{260}(3), 856--872
  (2017)

\bibitem{przybylski2010recursive}
Przybylski, A., Gandibleux, X., Ehrgott, M.: A recursive algorithm for finding
  all nondominated extreme points in the outcome set of a multiobjective
  integer programme.
\newblock INFORMS Journal on Computing \textbf{22}(3), 371--386 (2010)

\bibitem{Przybylski2010TwoPhase}
Przybylski, A., Gandibleux, X., Ehrgott, M.: A two phase method for
  multi-objective integer programming and its application to the assignment
  problem with three objectives.
\newblock Discrete Optimization \textbf{7}(3), 149--165 (2010)

\bibitem{Rodder1977GeneralizedSaddlepoint}
R{\"o}dder, W.: A generalized saddlepoint theory: Its application to duality
  theory for linear vector optimum problems.
\newblock European Journal of Operational Research \textbf{1}(1), 55--59 (1977)

\bibitem{sourd2008multiobjective}
Sourd, F., Spanjaard, O.: A multiobjective branch-and-bound framework:
  Application to the biobjective spanning tree problem.
\newblock INFORMS Journal on Computing \textbf{20}(3), 472--484 (2008)

\bibitem{Teghem2009}
Teghem, J.: Multi-objective integer linear programming.
\newblock In: C.A. Floudas, P.M. Pardalos (eds.) Encyclopedia of Optimization,
  pp. 2448--2454. Springer US, Boston, MA (2009)

\bibitem{ulungu1995two}
Ulungu, E.L., Teghem, J.: The two phases method: An efficient procedure to
  solve bi-objective combinatorial optimization problems.
\newblock Foundations of Computing and Decision Sciences \textbf{20}(2),
  149--165 (1995)

\bibitem{VINCENT2013498}
Vincent, T., Seipp, F., Ruzika, S., Przybylski, A., Gandibleux, X.: Multiple
  objective branch and bound for mixed 0-1 linear programming: Corrections and
  improvements for the biobjective case.
\newblock Computers \& Operations Research \textbf{40}(1), 498--509 (2013)

\bibitem{wolsey1981integer}
Wolsey, L.A.: Integer programming duality: Price functions and sensitivity
  analysis.
\newblock Mathematical Programming \textbf{20}(1), 173--195 (1981)

\bibitem{wolsey2014integer}
Wolsey, L.A., Nemhauser, G.L.: Integer and Combinatorial Optimization.
\newblock John Wiley \& Sons (2014)

\end{thebibliography}

\begin{appendix}
\normalsize
\section*{Appendix}
\section{Proofs of Results in Section \ref{sec:setorder}}
\label{sec:app_setorder}

The set-ordering in Definition \ref{def:Set_Order} is well-defined for subsets of $\reals^k$.
We extend it to $\pm M_{\infty}$ by defining $-M_{\infty} \preceqq S \preceqq M_{\infty}$ for all nonempty sets $S \subseteq \reals^k$, $-M_{\infty} \preceqq -M_{\infty}$, and $M_{\infty} \preceqq M_{\infty}$.
We further assume that $S \npreceqq -M_{\infty}$, and $M_{\infty} \npreceqq S$ for any $S \subseteq \reals^k$.
The relation ``$\preceqq$'' is thus defined on the extended power set of $\reals^k$ (excluding the empty set) and is transitive thereon, but neither reflexive nor antisymmetric. However, the relation defines a partial order on the family of sets $\mathcal{E}$ defined in \eqref{eq:setE}.

\PropMaxSinE*
\proof 
We prove the result for $\Max(S)$; the proof for $\Min(S)$ is similar and therefore omitted.
Suppose $S$ is nonempty and has points that are nondominated from above. 
Let $s,t \in \Max(S)$ with $s \neq t$. By definition of nondominance, $s \not \leq t$ and $t \not \leq s$. Thus, $\Max(A) \in \mathcal{E}$.
\qed
\endproof


\PropSetPartialOrder*
\proof 
We first show that the relations is reflexive. Consider a set $S \in \mathcal{E}$. If $S = \pm M_{\infty}$, then $S \preceqq S$ by definition of $\pm M_{\infty}$. 
Otherwise, if $s \in S$, then $s \leqq s$ by the reflexivity of $\leqq$. 
Moreover, if $t \in S$ with $s \not = t$, then $t \not \leq s$ because distinct elements of $S$ are incomparable as $S \in \mathcal{E}$. Therefore $S \preceqq S$.

Next, to see that $\preceqq$ is antisymmetric, consider $S,T \in \mathcal{E}$ with $S \preceqq T$ and $T \preceqq S$. If one of $S,T = \pm M_{\infty}$, then $S = T$. Otherwise, let $S,T \subseteq \reals^k$ and $s \in S$. Because $S \preceqq T$, there exists $t \in T$ such that $s \leqq t$. On the other hand, because $T \preceqq S$, there is $s' \in S$ such that $t \leqq s'$. Therefore, $s \leqq t \leqq s'$. 
But distinct elements of $S$ are incomparable, which implies that $s = s'$, so that $s = s' = t$. Therefore $S \subseteq T$. Similarly, $T \subseteq S$ as well, so that $S = T$.

Finally, to establish transitivity, suppose $S,T,U \in \mathcal{E}$ with $S \preceqq T$ and $T \preceqq U$. If $S = -M_{\infty}$ or $U = M_{\infty}$, the results holds trivially. If $S$ or $T = M_{\infty}$, then $U= M_{\infty}$ and the result holds. Similarly, if $T$ or $U$ equals $-M_{\infty}$, then so does $S$ and the result follows. Assume, therefore, that $S,T,U \neq \pm M_{\infty}$.

Let $s \in S$. Then, there exists $t \in T$ such that $s \leqq t$ and $u \in U$ such that $t \leqq u$. Therefore $s \leqq u$. 
On the other hand, suppose there exists $u \in U$ and $s \in S$ such that $u \leq s$. Then, there is $t \in T$ such that $s \leqq t$, which implies that $u \leq s \leqq t$. This contradicts the hypothesis that $T \preceqq U$. So, $s$ and $u$ must be incomparable. Therefore, $S \preceqq U$. 

Thus, the relation $\preceqq$ is reflexive, antisymmetric and transitive on $\mathcal{E}$ and therefore defines a partial order thereon.
\qed
\endproof

\LemOrderSubset*
\proof
If $U = M_{\infty}$, the result holds by definition. Therefore, suppose $U \neq M_{\infty}$, and let $t \in T$. Because $t \in S$, there is an element $u \in U$ such that $t \leqq u$.
Moreover, given $u \in U$, there is no $t \in T$ such that $u \leq t$ because that would contradict $S \preceqq U$. Thus, $T \preceqq U$.
\qed
\endproof

\LemOrderMax*
\proof
If $S = \emptyset$, then $\Max(S) = -M_{\infty} \preceqq U$. 
If $S$ is unbounded above, then we must have $U = M_{\infty}$ and $\Max(S) = M_{\infty} \preceqq U$. Suppose, therefore, that $S$ is nonempty and bounded above. The result then follows from Lemma \ref{lem:order_subset} by choosing $T = \Max(S)$.
\qed \endproof

\section{Omitted Proofs from Section \ref{sec:relaxation}} \label{sec:app_relaxation}

\PropMOLPRelax*
\proof 
If $\eqref{eq:MOIP}$ is infeasible, then $\Max(\mc{Y}_{\mathrm{MOIP}}) = -M_{\infty}$ and the result is trivially true. If \eqref{eq:MOIP} is unbounded then so is \eqref{eq:MOLP} because $\mathcal{Y}_{\mathrm{MOIP}}\subseteq \mathcal{Y}_{\mathrm{MOLP}}$. Finally, if $Cx^* \in \Max(\mc{Y}_{\mathrm{MOIP}})$, then $Cx^*$ is a feasible objective to \eqref{eq:MOLP}. Therefore either $\Max(\mc{Y}_{\mathrm{MOLP}}) = M_{\infty}$ or there exists $y \in \Max(\mc{Y}_{\mathrm{MOLP}})$ such that $Cx^* \leqq y$.

Now suppose $C\tilde{x} \in \Max(\mc{Y}_{\mathrm{MOLP}})$ such that $Cx^* \not = C\tilde{x}$. Suppose if possible that $C\tilde{x} \leq Cx^*$. Then, because $x^* \in \mathcal{X}$, it is also feasible to \eqref{eq:MOLP}, which contradicts the nondominance of $\tilde{x}$. Thus, $C\tilde{x} \not \leq Cx^*$ so that $Cx^*$ and $C\tilde{x}$ are incomparable. 
\qed
\endproof

\PropMOIPLeqqCH*

\proof The proof is similar to that of Proposition \ref{prop:MOLP_relax} and is omitted.  \endproof

\PropCHMOIPIntegral*

\proof
If $x^*$ is efficient for \eqref{eq:MOIP}, then it must be feasible to \eqref{eq:MOIP} and therefore integral. To see the opposite containment, suppose $x^*$ is efficient for \eqref{eq:CH} and integral. Then, $x^*$ is an integral point in the feasible region of \eqref{eq:CH} and therefore $x^{*} \in \mathcal{X}_{\mathrm{MOIP}}$. Suppose if possible that $x^*$ is not efficient for \eqref{eq:MOIP}.
Proposition \ref{prop:MOIP_leqq_CH} implies that for every feasible $x$ to \eqref{eq:MOIP}, there is a $y$ feasible to \eqref{eq:CH} such that $Cx \leqq Cy$. In particular, if $x^*$ is not efficient for \eqref{eq:MOIP}, then there is a $y$ feasible to \eqref{eq:CH} such that $Cx^*\leq Cy$. However, this contradicts the hypothesis that $x^{*}$ was efficient to \eqref{eq:CH}. Thus, $x^*$ must be efficient for \eqref{eq:MOIP}. 
\qed\endproof

\PropCHIntSol*
\proof Suppose that $\text{conv}(\mathcal{X}_{\mathrm{MOIP}})$ has a vertex. If an MOLP has efficient solutions, then at least one must occur at a vertex of the feasible region \cite[Theorem 4.3.8(ii)]{luc_2016_Book}. In particular, \eqref{eq:CH} has at least one efficient vertex. Because the vertices of the feasible region of \eqref{eq:CH} are all integral, this implies that there is at least one integral efficient solution to \eqref{eq:CH}. 

Now suppose that $\text{conv}(\mathcal{X}_{\mathrm{MOIP}})$ does not have vertices. Then, by \cite[Theorem 4.3.8(iii)]{luc_2016_Book}, it has a nonempty face $F$ such that every $x \in F$ is an efficient solution. Because $F$ is a face of $\text{conv}(\mathcal{X}_{\mathrm{MOIP}})$, there is a matrix $B$ and a vector $d$ such that $F = \{x \in \text{conv}(\mathcal{X}_{\mathrm{MOIP}})\ \vert \ Bx = d\}$ and such that $Bx \leqq d$ for all $x \in \mathcal{X}_{\mathrm{MOIP}}$. If $x^* \in F$, then $x^* = \sum_{i = 1}^{p} t^ix^i$ for some positive integer $p$, $t^1, t^2, \ldots,t^p \in [0,1]$ and $x^1,x^2, \ldots, x^p\in \mathcal{X}_{\mathrm{MOIP}}$ with $\sum_{i = 1}^{p}t^i = 1$ because $F \subseteq \text{conv}(\mathcal{X}_{\mathrm{MOIP}})$. Then, 
\[
Bx^* = B\Big(\sum_{i = 1}^{p}t^ix^i\Big)=\sum_{i = 1}^{p}t^iBx^i = d.
\]
Because $Bx^i\leqq d$ for each $i$, this implies that $Bx^i = d$ for each $i$. So, $x^i\in F$ and therefore \eqref{eq:CH} has an integral efficient solution. 
\qed\endproof

\PropScalConvex*

\proof
Let $x^*$ be an efficient solution of \eqref{eq:CH}. By Lemma \ref{lem:MOLP_scal_iff}, there is a scalarizing vector $\mu\in \reals_{>}$ such that $x^*$ is optimal to 
\begin{equation}
\label{eq:scal_CH}\tag{CH$_\mu$}
\begin{aligned}
\max & \enskip \mu^\top Cx\\
\text{s.t. } &\enskip x \in \text{conv}(\{x\in \Z_+^n\ \vert \ Ax \leqq b\}).
\end{aligned}
\end{equation}
\noindent On the other hand, because $x^*$ is feasible to \eqref{eq:MOIP}, it is feasible to the IP
\begin{equation}
\label{eq:scal_MOIP}\tag{MOIP$_\mu$}
\begin{aligned}
\max & \enskip \mu^\top Cx\\
\text{s.t. } & \enskip x \in \{ x \in \Z_+^n \ \vert \ Ax \leqq b \}.
\end{aligned}
\end{equation}

Note that \eqref{eq:scal_CH} is the convex-hull relaxation of the single-objective IP \eqref{eq:scal_MOIP}. Thus, $x^*$ must be optimal to \eqref{eq:scal_MOIP}. Therefore, $x^*$ is a supported efficient solution to \eqref{eq:MOIP}.

Conversely, suppose $x^*$ is a supported efficient solution of \eqref{eq:MOIP}. Then, there exists a scalarizing vector $\mu\in \reals^k_>$ such that $x^*$ is an optimal solution to \eqref{eq:scal_MOIP}. Once again, because \eqref{eq:scal_CH} is the convex hull relaxation of \eqref{eq:scal_MOIP}, $x^*$ is optimal to \eqref{eq:scal_CH} as well. It follows from Lemma \ref{lem:MOLP_scal_iff} that $x^*$ is efficient for \eqref{eq:CH}.
\qed
\endproof




\section{An Example to Illustrate Theorem \ref{thm:MO_Lag_iff_1}} \label{sec:ex:LD_not_strong_supp}

This example illustrates Theorem \ref{thm:MO_Lag_iff_1}. We consider an MOIP whose feasible region does not satisfy condition \eqref{eq:FR_LAG} and show that the Lagrangian dual for this problem is not strong at a supported solution.

\begin{example}\label{ex:LD_not_strong_supp}
Consider the problem
\begin{equation*} \begin{aligned}
\begin{aligned}
\max \enskip & \begin{bmatrix} 1 & 0 \\ 0 & 1\end{bmatrix}\begin{bmatrix} x_1\\x_2 \end{bmatrix} \quad 
\text{s.t. } \enskip 2x_1 + 4x_2 \leqq 5,\  4x_1 + 2x_2 \leqq 5,\  x_1, x_2 \in \{0,1\}.
\end{aligned}
\end{aligned} \end{equation*}

For this problem, $\mathcal{X} = \mathcal{Y} = \{(0,0)^\top,(1,0)^\top,(0,1)^\top\}$ and the supported nondominated points are $(1,0)^\top$ and $(0,1)^\top$. Set $Q = \{0,1\}^2$, $A^1 = \begin{bmatrix} 2 & 4\\ 4 & 2 \end{bmatrix}$, and $b^1 = (5, 5)^\top$. Then,
\begin{equation*} \begin{aligned}
\text{conv}(Q \cap \{x \in \mathbb{R}^n\ \vert \ A^1x \leqq b^1\}) \subset \text{conv}(Q) \cap \{x \in \mathbb{R}^n\ \vert \ A^1x \leqq b^1\},
\end{aligned} \end{equation*}
\noindent where the containment is strict. Enumerating $x \in Q$, we have the Lagrangian dual
\begin{equation*} \begin{aligned}
\Min(\mathcal{Y}_{\mathrm{LD}}) = \Min \Bigg( \bigcup\limits_{\Lambda \geqq 0} \Max \Bigg\{
& \begin{pmatrix} 5\lambda_{11} + 5\lambda_{12}\\ 5\lambda_{21} + 5\lambda_{22} \end{pmatrix},
\begin{pmatrix} 3\lambda_{11} + \lambda_{12}\\ 1 + \lambda_{21} + 3\lambda_{22} \end{pmatrix}, \\
& 
\begin{pmatrix} 1 +\lambda_{11} + 3\lambda_{12}\\ 3\lambda_{21} + \lambda_{22} \end{pmatrix},
\begin{pmatrix} 1-\lambda_{11} -\lambda_{12}\\ 1-\lambda_{21} -\lambda_{22} \end{pmatrix}
\Bigg\} \Bigg).
\end{aligned} \end{equation*}

Then, there is no $\Lambda \in \reals^{2\times 2}_{+}$ such that $(1,0)^\top \in \Min(\mathcal{Y}_{\mathrm{LD}})$. 
To see this, first suppose that $\Lambda$ was such that $\begin{pmatrix} 5\lambda_{11} + 5\lambda_{12}\\ 5\lambda_{21} + 5\lambda_{22} \end{pmatrix} = (1,0)^\top$. Then, we must have $\lambda_{11} + \lambda_{12} = \frac{1}{5}$ and $\lambda_{21} = \lambda_{22} = 0$. But then, because one of $\lambda_{11},\lambda_{12}$ must be positive, there are positive parameters $\delta_1 = 3\lambda_{11} + \lambda_{12}$ and $\delta_2 = \lambda_{11} + 3\lambda_{22}$ with $\delta_1,\delta_2 \geq \frac{1}{5}$ such that the relaxation \eqref{eq:Lagrangian_Relaxation} is 
\begin{equation*} 
\begin{aligned}
\Max \left\{\begin{pmatrix}1\\0 \end{pmatrix},
\begin{pmatrix}\delta_1\\1 \end{pmatrix},
\begin{pmatrix}1+\delta_2\\0 \end{pmatrix},
\begin{pmatrix}\frac{4}{5}\\1 \end{pmatrix}
\right\}.
\end{aligned} 
\end{equation*}
\noindent 
Then, $1+\delta_2 \geq \frac{6}{5} > 1$, which implies that $(1,0)\not \in \Max(\mathcal{Y}_{\mathrm{LR}(\Lambda)})$. 

Because $\lambda_{21},\lambda_{22} \geq 0$, there is no feasible $\Lambda$ such that $1+ \lambda_{21} +3\lambda_{22} = 0$. So, $\begin{pmatrix} 3\lambda_{11} + \lambda_{12}\\ 1 + \lambda_{21} + 3\lambda_{22} \end{pmatrix} \not = \begin{pmatrix}1\\0 \end{pmatrix}.$

Next, if $\begin{pmatrix} 1 +\lambda_{11} + 3\lambda_{12}\\ 3\lambda_{21} + \lambda_{22} \end{pmatrix} = \begin{pmatrix}1\\0 \end{pmatrix}$, then $\Lambda = 0$. However, $(1,0) \not \in \Max(\mathcal{Y}_{\mathrm{LR}(0)})$ because
\begin{equation*} 
\begin{aligned}\Max(\mathcal{Y}_{\mathrm{LR}(0)})=\Max\{(0,0),(0,1),(1,0),(1,1)
\} = \{(1,1)\}.
\end{aligned} 
\end{equation*}

Finally, if $\Lambda$ was such that $\begin{pmatrix} 1-\lambda_{11} -\lambda_{12}\\ 1-\lambda_{21} -\lambda_{22} \end{pmatrix} = \begin{pmatrix} 1\\0\end{pmatrix}$ then $\lambda_{11} = \lambda_{12} = 0$ and $\lambda_{21} + \lambda_{22} = 1$. Then, there are positive parameters $\delta_1 = \lambda_{21} + 3\lambda_{22}$ and $\delta_2 = 3\lambda_{21} + \lambda_{22}$ such that $\delta_1, \delta_2 \geq 1$ and 
\begin{equation*} \begin{aligned}
\Max(\mathcal{Y}_{\mathrm{LR}(\Lambda)}) = \Max\left\{
\begin{pmatrix} 0\\ 5 \end{pmatrix},
\begin{pmatrix} 0\\ 1+\delta_1 \end{pmatrix},
\begin{pmatrix}1 \\ \delta_2 \end{pmatrix}, 
\begin{pmatrix} 1\\ 0 \end{pmatrix}
\right\},
\end{aligned} 
\end{equation*}
\noindent so that $(1,0) \not \in \Max(\mathcal{Y}_{\mathrm{LR}(\Lambda)})$ because $\delta_2 \geq 1$.
Therefore, there are no Lagrangian relaxations which are tight at $(1,0)^\top$. Moreover, the relaxations are bounded away from $(1,0)^\top$ so that $(1,0)^\top$ cannot be a limit point of $\mathcal{Y}_{\mathrm{LD}}$.

Thus, condition \eqref{eq:FR_LAG} is not satisfied, and the Lagrangian dual for this problem is not strong at a supported efficient solution.
\qed
\end{example}

\end{appendix}

\end{document}